\documentclass{article}
\usepackage{a4,amsmath, amssymb, amscd}

\DeclareMathSymbol{\Z}{\mathalpha}{AMSb}{"5A}
\DeclareMathSymbol{\PP}{\mathalpha}{AMSb}{"50}
\DeclareMathSymbol{\Q}{\mathalpha}{AMSb}{"51}
\DeclareMathSymbol{\N}{\mathalpha}{AMSb}{"4E}
\DeclareMathSymbol{\R}{\mathalpha}{AMSb}{"52}

\newcommand {\ov} {\overline}

\newcommand {\Hom} {{\rm Hom}}
\newcommand {\im} {{\rm Im}\,}
\newcommand {\Ker} {{\rm Ker}\,}

\newcommand {\Rad} {\mbox{Rad\,}}
\newcommand {\Pic} {{\rm Pic}\,}
\newcommand {\Out} {{\rm Out\,}}
\newcommand {\Aut} {{\rm Aut\,}}
\newcommand {\Inn} {{\rm Inn\,}}
\newcommand {\car} {{\rm char\,}}
\newcommand {\ie} {{\em i.e.}\ }
\newcommand {\Ext} {{\rm Ext\,}}
\newcommand {\Der} {{\rm Der\,}}
\newcommand {\Ad} {{\rm Ad\,}}
\newcommand {\id} {{\rm id}\,}
\newcommand {\ad} {{\rm ad}\,}

\newenvironment {Bew}%
       {\noindent {\it Proof:}}
       {\begin{flushright}{\it q.e.d.} \theDef\end{flushright}}
\newtheorem {Def} {Definition} [section]

\newtheorem {Prop} [Def] {Proposition} 

\newtheorem {Kor} [Def] {Corollary}
\newtheorem {Thm} [Def] {Theorem}
\newtheorem {Lem} [Def] {Lemma}
\newtheorem {Bem} [Def] {Remark}
\newtheorem {Beme} [Def] {Remarks}
\newtheorem {Not} [Def] {Notations}
\date{}
\begin{document}

\title{The Lie algebra structure on the first Hochschild cohomology group of a monomial algebra}
\author{Claudia Strametz
        \\
        \\ \small\em D\'epartement de Math\'ematiques, Universit\'e de Montpellier 2,
        \\ \small\em F-34095 Montpellier cedex 5, France
        \\ \small E-mail: strametz@math.univ-montp2.fr}
\maketitle
\renewcommand\abstractname{R\'esum\'e}
\begin{abstract}
 Nous \'etudions la structure d'alg\`{e}bre de Lie du premier groupe de la 
 cohomologie de Hochschild d'une alg\`{e}bre monomiale de dimension finie
 $\Lambda$, en termes combinatoires de son carquois, en quelconque 
 caract\'eristique. Cela nous permet aussi d'examiner la composante de 
 l'identit\'e du groupe alg\'ebrique des automorphismes ext\'erieurs de
 $\Lambda$ en caract\'eristique z\'ero. Nous donnons des crit\`{e}res pour la 
 (semi-)simplicit\' et la r\'esolubilit\'e.
\end{abstract}
\renewcommand\abstractname{Abstract}
\begin{abstract}
 We study the Lie algebra structure of the first Hochschild cohomology group 
 of a finite dimensional monomial algebra $\Lambda$, in terms of the 
 combinatorics of its quiver, in any characteristic. This allows us also to 
 examine the identity component of the algebraic group of outer automorphisms 
 of $\Lambda$ in characteristic zero. Criteria for the (semi-)simplicity, the
 solvability, the reductivity, the commutativity and the nilpotency are given.
\end{abstract}
\bigskip
{\bf 2000 Mathematics Subject Classification:} 16E40 16W20

\bigskip
{\bf Keywords}: cohomology, Hochschild, Lie, monomial

\section{Introduction}

The Hochschild cohomology ${\rm H}^*(\Lambda,\Lambda)$ of any associative 
algebra $\Lambda$ over a field $k$ has the structure of a Gerstenhaber 
algebra (see \cite{Ge}). In particular, the first Hochschild cohomology 
group ${\rm H}^1(\Lambda,\Lambda)$ is a Lie algebra, a fact which can be 
verified 
directly. Note that in the finite dimensional case in characteristic $0$ this 
Lie algebra can be regarded as the Lie algebra of the algebraic group of outer
automorphisms 
$\Out(\Lambda)=\Aut(\Lambda)/\Inn(\Lambda)$ of the algebra $\Lambda$. 
It has been treated by Guil-Asensio and Saor\'{\i}n in \cite{GAS1}. 
Huisgen-Zimmermann and Saor\'{\i}n proved in \cite {HZS} that the identity 
component $\Out(\Lambda)^\circ$ of the outer automorphism group of $\Lambda$ is
invariant under derived equivalence and Keller \cite{Ke} showed that 
Hochschild cohomology is preserved under 
derived equivalence as a graded (super) Lie algebra. Consequently the Lie algebra 
${\rm H}^1(\Lambda,\Lambda)$ is invariant under derived equivalence. 

The purpose of this paper is to study the Lie algebra 
structure of ${\rm H}^1(\Lambda,\Lambda)$ in the case of finite
dimensional monomial algebras without any restriction on the characteristic of
the field $k$ using only algebraic methods. The relationship between 
${\rm H}^1(\Lambda,\Lambda)$ and $\Out(\Lambda)$ allows us to transfer the 
results obtained in this way to the identity component $\Out(\Lambda)^\circ$
of the algebraic group $\Out(\Lambda)$ in characteristic $0$. 
Thus we give a different proof of Guil-Asensio and Saor\'{\i}n's criterion for the
solvability of $\Out(\Lambda)^\circ$ and we generalize some results they 
obtained using algebraic group theory and methods in algebraic geometry.

This paper is organized in the following way: in the first section we will use
the minimal projective resolution of a monomial algebra $\Lambda$ (as a 
$\Lambda$-bimodule) given by Bardzell in \cite{Ba} to get a handy description 
of the Lie algebra ${\rm H}^1(\Lambda,\Lambda)$ in terms of parallel paths. 
The purpose of the second section is to link this description to Guil-Asensio
and Saor\'{\i}n's work on the algebraic group of outer automorphisms of monomial
algebras (see \cite{GAS2}). 
In section three we carry out the study of the Lie algebra 
${\rm H}^1(\Lambda,\Lambda)$. In particular we 
give criteria for the (semi-)simplicity, the solvability, the reductivity, 
the commutativity and the nilpotency of this Lie algebra and consequently of 
the connected algebraic group $\Out(\Lambda)^\circ$, in terms of the 
combinatorics of the quiver of $\Lambda$. Finally a Morita equivalence given by
Gabriel in \cite{Ga} in case $\car k=p>0$ between group algebras $kG$ where 
the finite group $G$ admits a normal cyclic Sylow $p$-subgroup and certain 
monomial algebras allows us to apply some of our results to the Lie algebra 
${\rm H}^1(kG,kG)$. Note that the Hochschild cohomology is Morita invariant as
a Gerstenhaber algebra (see \cite{GS}).

We give some notation and terminology which we keep throughout the paper. Let 
$Q$ denote a finite quiver (that is a finite oriented graph) and $k$ an 
algebraically closed field. For all $n\in\N$ let $Q_n$ be the set of oriented 
paths of length $n$ of $Q$. Note that $Q_0$ is the set of vertices and that 
$Q_1$ is the set of arrows of $Q$. We denote by $s(\gamma)$ the source vertex 
of an (oriented) path $\gamma$ of $Q$ and by $t(\gamma)$ its terminus vertex. 
The path algebra $kQ$ is the $k$-linear span of the set of paths of $Q$ where 
the multiplication of $\beta\in Q_i$ and $\alpha\in Q_j$ is provided by the 
concatenation $\beta\alpha\in Q_{i+j}$ if $t(\alpha)=s(\beta)$ and $0$ otherwise.

We denote by $\Lambda$ a finite dimensional monomial $k$-algebra, that is a 
finite dimensional $k$-algebra which is 
isomorphic to a quotient of a path algebra $kQ/I$ where the two-sided ideal 
$I$ of $kQ$ is generated by a set $Z$ of paths of length $\ge2$. We shall 
assume that $Z$ is minimal, \ie no proper subpath of a path in $Z$ is again in
$Z$. Let $B$ be the set of paths of $Q$ which do not contain any element of $Z$
as a subpath. It is clear that the (classes modulo $I=\langle Z\rangle$ of)
elements of $B$ form a basis of $\Lambda$. We shall denote by $B_n$ the subset
$Q_n\cap B$ of $B$ formed by the paths of length $n$.

Let $E\simeq kQ_0$ be the separable subalgebra of $\Lambda$ generated by the
(classes modulo $I$ of the) vertices of $Q$. We have a Wedderburn-Malcev 
decomposition $\Lambda=E\oplus r$ where $r$ denotes the Jacobson radical of 
$\Lambda$.

This work will form part of a Ph.D. thesis under the supervision of Claude 
Cibils. I would like to thank him for his comments and his encouragement which
I appreciated.

\section{Projective resolutions and the Lie bracket}

The Hochschild cohomology 
${\rm H}^*(\Lambda,\Lambda)=\Ext_{\Lambda^e}^*(\Lambda,\Lambda)$ of a $k$-algebra $\Lambda$ can be computed using different projective resolutions of 
$\Lambda$ over its enveloping algebra 
$\Lambda^e=\Lambda\otimes_k\Lambda^{\rm op}$. The standard resolution 
$\mathcal{P}_{\rm Hoch}$ is
$$\cdots\rightarrow\Lambda^{\otimes^n_k}\stackrel{\delta}{\rightarrow}
   \Lambda^{\otimes^{n-1}_k}\rightarrow\cdots\rightarrow
   \Lambda\otimes_k\Lambda\stackrel{\varepsilon}{\rightarrow}\Lambda
   \rightarrow0$$
where $\varepsilon(x_1\otimes_kx_2)=x_1x_2$ and
$$\delta(x_1\otimes_k\cdots\otimes_kx_n)=\sum_{i=1}^{n-1}(-1)^{i+1}
   x_1\otimes_k\cdots\otimes_kx_ix_{i+1}\otimes_k\cdots\otimes_kx_n$$
for $x_1,\ldots,x_n\in\Lambda$. Applying the functor 
$\Hom_{\Lambda^e}(\_\,,\Lambda)$ to $\mathcal{P}_{\rm Hoch}$ and
identifying 
$\Hom_{\Lambda^e}(\Lambda\otimes_k\Lambda^{\otimes_k^n}\otimes_k\Lambda,\Lambda)$ 
with
$\Hom_k(\Lambda^{\otimes_k^n},\Lambda)$ for all $n\in\N$, yields the cochain complex 
$\mathcal{C}_{\rm Hoch}$ defined by Hochschild:
$$0\rightarrow\Lambda\stackrel{d_0}{\longrightarrow}\Hom_k(\Lambda,\Lambda)
  \rightarrow\cdots\rightarrow\Hom_k(\Lambda^{\otimes_k^n},\Lambda)
  \stackrel{d_n}{\longrightarrow}\Hom_k(\Lambda^{\otimes_k^{n+1}},\Lambda)
  \rightarrow\cdots$$
where $(d_0a)(x)=ax-xa$ for all $a,x\in\Lambda$ and
$$\begin{array}{rcl}
  (d_nf)(x_1\otimes_k\cdots\otimes_kx_{n+1})&=&x_1f(x_2\otimes_k\cdots\otimes_k
  x_{n+1})+\\
  &&\sum_{i=1}^n(-1)^if(x_1\otimes_k\cdots\otimes_kx_ix_{i+1}\otimes_k\cdots
  \otimes_kx_{n+1})\\
  &&+(-1)^{n+1}f(x_1\otimes_k\cdots\otimes_kx_n)x_{n+1}
\end{array}$$
for all $f\in\Hom_k(\Lambda^{\otimes_k^n},\Lambda)$, $n\in\N$, and
$x_1,\ldots,x_{n+1}\in\Lambda$. In 1962 the structure of a Gerstenhaber 
algebra was introduced on the Hochschild cohomology 
${\rm H}^*(\Lambda,\Lambda)$ by Gerstenhaber in \cite{Ge}. In particular the 
first cohomology group ${\rm H}^1(\Lambda,\Lambda)$ which is the quotient of 
the derivations modulo the inner derivations of $\Lambda$ is a Lie algebra 
whose bracket is induced by the Lie bracket on $\Hom_k(\Lambda,\Lambda)$
$$[f,g]=f\circ g-g\circ f$$
where $f,g\in\Hom_k(\Lambda,\Lambda)$. The $1$-coboundaries, \ie the 
derivations $\Der_k(\Lambda)=\Ker d_1=\{f\in\Hom_k(\Lambda,\Lambda)\mid 
f(ab)=af(b)+f(a)b\ \forall a,b\in\Lambda\}$ form a Lie subalgebra of 
$\Hom_k(\Lambda,\Lambda)$ and the $1$-cocycles, \ie the inner derivations 
$\Ad_k(\Lambda)=\im d_0=\{f\in\Hom_k(\Lambda,\Lambda)\mid\exists a\in\Lambda
\ {\rm such\ that}\ f(x)=ax-xa\ \forall x\in\Lambda\}$ an ideal of this Lie algebra.

In order to compute the first Hochschild cohomology group 
${\rm H}^1(\Lambda,\Lambda)$ of the finite dimensional monomial algebra 
$\Lambda$, we shall use the minimal projective resolution of the 
$\Lambda$-bimodule $\Lambda$ given by Bardzell in \cite{Ba}. The part of this 
resolution $\mathcal{P}_{\rm min}$ in which we are interested is the following:
$$\cdots\rightarrow\Lambda\otimes_EkZ\otimes_E\Lambda\stackrel{\delta_1}
  {\longrightarrow}\Lambda\otimes_EkQ_1\otimes_E\Lambda\stackrel{\delta_0}
  {\longrightarrow}\Lambda\otimes_E\Lambda\stackrel{\pi}{\longrightarrow}
  \Lambda\rightarrow0$$
where the $\Lambda$-bimodule morphisms are given by
$$\begin{array}{rcl}
 \pi(\lambda\otimes_E\mu)&=&\lambda\mu\\
 \delta_0(\lambda\otimes_Ea\otimes_E\mu)&=&\lambda a\otimes_E\mu-
   \lambda\otimes_Ea\mu\qquad{\rm and}\\
 \delta_1(\lambda\otimes_Ep\otimes_E\mu)&=&\sum_{d=1}^n\lambda a_n\ldots 
 a_{d+1}\otimes_Ea_d\otimes_Ea_{d-1}\ldots a_1\mu
\end{array}$$
for all $\lambda,\mu\in\Lambda$, $a,a_n,\ldots,a_1\in Q_1$ and 
$p=a_n\ldots a_1\in Z$ (with the conventions $a_{n+1}=t(a_n)$ and 
$a_0=s(a_1))$.

\begin{Bem}
Our description of $\mathcal{P}_{\rm min}$ is equivalent to Bardzell's, 
because if $X$ denotes a set of paths and if $kX$ is the corresponding 
$E$-bimodule, then the map 
$\oplus_{\gamma\in X}\Lambda t(\gamma)\otimes_ks(\gamma)\Lambda\rightarrow
\Lambda\otimes_EkX\otimes_E\Lambda$ which is given by 
$\lambda t(\gamma)\otimes_k s(\gamma)\mu\mapsto
\lambda\otimes_E\gamma\otimes_E\mu$, where $\gamma\in X$ and 
$\lambda,\mu\in\Lambda$, is clearly a $\Lambda$-bimodule isomorphism. 
Note
that the $\Lambda$-bimodules $\Lambda\otimes_E\Lambda$ and 
$\Lambda\otimes_EE\otimes_E\Lambda\simeq\Lambda\otimes_EkQ_0\otimes_E\Lambda$
are isomorphic.
\end{Bem}

\begin{Lem}
Let $M$ be an $E$-bimodule and $T$ a $\Lambda$-bimodule. Then the vector space
$\Hom_{\Lambda^e}(\Lambda\otimes_EM\otimes_E\Lambda, T)$ is isomorphic to
$\Hom_{E^e}(M,T)$.
\end{Lem}

\begin{Bew}
 The linear morphism $\Hom_{\Lambda^e}(\Lambda\otimes_EM\otimes_E\Lambda,T)
 \rightarrow\Hom_{E^e}(M,T)$ which sends a $\Lambda$-bimodule morphism
 $f:\Lambda\otimes_EM\otimes_E\Lambda\rightarrow T$ to the $E$-bimodule 
 morphism $M\rightarrow T$ given by 
 $m\mapsto f(1_\Lambda\otimes_Em\otimes_E1_\Lambda)$ and the linear morphism
 $\Hom_{E^e}(M,T)\rightarrow
 \Hom_{\Lambda^e}(\Lambda\otimes_EM\otimes_E\Lambda,T)$ which associates to
 $g:M\rightarrow T$ the element 
 $\Lambda\otimes_EM\otimes_E\Lambda\rightarrow T,\ 
 \lambda\otimes_Em\otimes_E\mu\mapsto\lambda g(m)\mu$, are inverse to each other.
\end{Bew}

Applying the functor $\Hom_{\Lambda^e}(\_\,,\Lambda)$ to 
$\mathcal{P}_{\rm min}$ and using the preceding Lemma yields the cochain 
complex $\mathcal{C}_{\rm min}$
$$0\longrightarrow\Hom_{E^e}(kQ_0,\Lambda)\stackrel{\delta^*_0}
  {\longrightarrow}
  \Hom_{E^e}(kQ_1,\Lambda)\stackrel{\delta^*_1}{\longrightarrow}
  \Hom_{E^e}(kZ,\Lambda)\longrightarrow\cdots$$
where the coboundaries $\delta^*_0$ and $\delta_1^*$ are given by
\begin{eqnarray*}
 (\delta^*_0f)(a)&=&af(s(a))-f(t(a))a\qquad\qquad\qquad{\rm and}\\
 (\delta^*_1g)(p)&=&\sum_{d=1}^na_n\ldots a_{d+1}g(a_d)a_{d-1}\ldots a_1
\end{eqnarray*}
where $f\in\Hom_{E^e}(kQ_0,\Lambda)$, $a,a_n,\ldots,a_1\in Q_1$, 
$g\in\Hom_{E^e}(kQ_1,\Lambda)$ and
$p=a_n\ldots a_1\in Z$.

As it is useful to interpret ${\rm H}^1(\Lambda,\Lambda)$ in terms of parallel
paths we introduce now the following notion: two paths $\varepsilon,\gamma$ of
$Q$ are called parallel if $s(\varepsilon)=s(\gamma)$ and 
$t(\varepsilon)=t(\gamma)$. If $X$ and $Y$ are sets of paths of 
$Q$, the set $X//Y$ of parallel paths is formed by the couples 
$(\varepsilon,\gamma)\in X\times Y$ such that $\varepsilon$ and 
$\gamma$ are parallel paths.
 For instance, $Q_0//Q_n$ is the set of oriented 
cycles of $Q$ of length $n$.

\begin{Lem}
 Let $X$ and $Y$ be sets of paths of $Q$ and let $kX$ and $kY$ be the 
 corresponding 
 $E$-bimodules. Then the vector spaces $k(X//Y)$ and $\Hom_{E^e}(kX,kY)$
 are isomorphic.
\end{Lem}

\begin{Bew}
 Define a linear morphism $k(X//Y)\rightarrow\Hom_{E^e}(kX,kY)$ by 
 sending
 $(\varepsilon,\gamma)\in X//Y$ to the elementary map which associates $\gamma$
 to $\varepsilon$ and $0$ to any other path of $X$. If $f:kX\rightarrow kY$ is
 an $E$-bimodule morphism, we have for every path  $\varepsilon\in X$ that
 $f(\varepsilon)=f(t(\varepsilon)\varepsilon s(\varepsilon))=t(\varepsilon)
 f(\varepsilon)s(\varepsilon)$ and thus
 $f(\varepsilon)=\sum_{(\varepsilon,\gamma)\in X//Y}
 \lambda_{\varepsilon,\gamma}\gamma$, where $\lambda_{\varepsilon,\gamma}\in k$. 
 This allows us to define a linear 
 morphism $\Hom_{E^e}(kX,kY)\rightarrow k(X//Y)$ which assigns to 
 $f:kX\rightarrow kY$ the element $\sum_{(\varepsilon,\gamma)\in X//Y}
 \lambda_{\varepsilon,\gamma}(\varepsilon,\gamma)$. Obviously, the two 
 morphisms are inverse to each other.
\end{Bew}

Using the following notations we will rewrite the coboundaries:

\begin{Not}
\ \begin{enumerate}
 \item Let $\varepsilon$ be a path in $Q$ and $(a,\gamma)\in Q_1//B$. We 
  denote by $\varepsilon^{(a,\gamma)}$ the sum of all nonzero paths 
  (\ie paths in $B$) obtained by replacing one appearance of the arrow $a$ in 
  $\varepsilon$ by the path $\gamma$. If the path $\varepsilon$ does not contain
  the arrow $a$ or if every replacement of $a$ in $\varepsilon$ by $\gamma$ is
  not a path in $B$, we
  set $\varepsilon^{(a,\gamma)}=0$. Suppose that $\varepsilon^{(a,\gamma)}=
  \sum_{i=1}^n\varepsilon_i$, where $\varepsilon_i\in B$ and let $\eta$ be a 
  path of $B$ parallel to $\varepsilon$. By abuse of language we denote by 
  $(\eta,\varepsilon^{(a,\gamma)})$ the sum $\sum_{i=1}^n(\eta,\varepsilon_i)$ 
  in $k(Z//B)$ (with the convention that $(\eta,\varepsilon^{(a,\gamma)})=0$ if
  $\varepsilon^{(a,\gamma)}=0$).
 \item The function $\chi_B:\prod_{n\in\N}Q_n\rightarrow\{0,1\}$ denotes the 
  indicator function which associates $1$ to each path $\gamma\in B$ and $0$ 
  to $\gamma\not\in B$.
 \item If $X$ is a set of paths of $Q$ and $e$ a vertex of $Q$, the set $Xe$ 
  is formed by the paths of $X$ with source vertex $e$. In the same way $eX$ 
  denotes the set of all paths of $X$ with terminus vertex $e$.
\end{enumerate}
\end{Not}

\begin{Bem}
\label{distinct}
 \ \\ If $(a,\gamma)\in Q_1//(B-Q_0)$, then all the nonzero summands of 
 $\varepsilon^{(a,\gamma)}$ are distinct.
\end{Bem}

If we carry out the identification suggested in the preceding Lemma, we obtain:

\begin{Prop}
\ \\ The beginning of the cochain complex $\mathcal{C}_{\rm min}$ can be characterized 
 in the following way
 $$0\rightarrow k(Q_0//B)\stackrel{\psi_0}{\longrightarrow}k(Q_1//B)
   \stackrel{\psi_1}{\longrightarrow}k(Z//B)
   \stackrel{\psi_2}{\longrightarrow}\cdots$$
 where the maps are given by
 $$\begin{array}{rrcl}
   \psi_0:&k(Q_0//B)&\longrightarrow&k(Q_1//B)\\
   &(e,\gamma)&\longmapsto&\sum_{a\in Q_1e}\chi_B(a\gamma)(a,a\gamma)-
                           \sum_{a\in eQ_1}\chi_B(\gamma a)(a,\gamma a)\\
   &&&\\
   \psi_1:&k(Q_1//B)&\longrightarrow&k(Z//B)\\
   &(a,\gamma)&\longmapsto&\sum_{p\in Z}(p,p^{(a,\gamma)})
 \end{array}$$
 In particular, we have 
 ${\rm H}^1(\Lambda,\Lambda)\simeq\Ker\psi_1/\im\psi_0$.
\end{Prop}

The verifications are straightforward.

\begin{Thm}
\label{bracket}
 The bracket
 $$[(a,\gamma),(b,\varepsilon)]=(b,\varepsilon^{(a,\gamma)})-
                                 (a,\gamma^{(b,\varepsilon)})$$
 for all $(a,\gamma), (b,\varepsilon)\in Q_1//B$ induces a Lie algebra 
 structure on $\Ker\psi_1/\im\psi_0$ such that ${\rm H}^1(\Lambda,\Lambda)$ 
 and $\Ker\psi_1/\im\psi_0$ are isomorphic Lie algebras.
\end{Thm}

\begin{Bew}
 As $\mathcal{P}_{\rm Hoch}$ and $\mathcal{P}_{\rm min}$ are projective 
 resolutions of the $\Lambda$-bimodule $\Lambda$, there exist, thanks to the 
 Comparison Theorem, chain maps $\omega:\mathcal{P}_{\rm Hoch}\rightarrow
 \mathcal{P}_{\rm min}$ and 
 \linebreak
 $\varsigma:\mathcal{P}_{\rm min}\rightarrow
 \mathcal{P}_{\rm Hoch}$. Let us choose these $\Lambda$-bimodule morphisms 
 such that
 $$\begin{array}{rrcl}
   \varsigma_0:&\Lambda\otimes_E\Lambda&\rightarrow&\Lambda\otimes_k\Lambda\\
   &\lambda\otimes_E\mu&\mapsto&\lambda\otimes_k\mu\\
   \varsigma_1:&\Lambda\otimes_EkQ_1\otimes_E\Lambda&\rightarrow&
   \Lambda\otimes_k\Lambda\otimes_k\Lambda\\
   &\lambda\otimes_Ea\otimes_E\mu&\mapsto&\lambda\otimes_ka\otimes_k\mu\\
   \omega_0:&\Lambda\otimes_k\Lambda&\rightarrow&\Lambda\otimes_E\Lambda\\
   &\lambda\otimes_k\mu&\mapsto&\lambda\otimes_E\mu\\
   \omega_1:&\Lambda\otimes_k\Lambda\otimes_k\Lambda&\rightarrow&
   \Lambda\otimes_EkQ_1\otimes_E\Lambda\\
   &\lambda\otimes_ka_n\ldots a_1\otimes_k\mu&\mapsto&\sum_{d=1}^n
   \lambda a_n\ldots a_{d+1}\otimes_Ea_d\otimes_Ea_{d-1}\ldots a_1\mu
 \end{array}$$
 where $\lambda,\mu\in\Lambda$, $a,a_n,\ldots,a_1\in Q_1$ and 
 $a_n\ldots a_1\in B$. As a chain map 
 between projective resolutions is unique up to chain homotopy equivalence, 
 $\omega$ and $\varsigma$ are unique up to chain homotopy equivalence and such
 that $\omega\circ\varsigma$ is homotopic to $\id_{\mathcal{P}_{\rm min}}$ and
 $\varsigma\circ\omega$ is homotopic to $\id_{\mathcal{P}_{\rm Hoch}}$. 
 Therefore, the cochain maps $\Hom_{\Lambda^e}(\omega,\Lambda)$ and 
 $\Hom_{\Lambda^e}(\varsigma,\Lambda)$ are such that 
 $\Hom_{\Lambda^e}(\varsigma,\Lambda)\circ\Hom_{\Lambda^e}(\omega,\Lambda)$ is
 homotopic to $\id_{\mathcal{C}_{\rm min}}$ and 
 $\Hom_{\Lambda^e}(\omega,\Lambda)\circ\Hom_{\Lambda^e}(\varsigma,\Lambda)$ is
 homotopic to $\id_{\Hom_{\Lambda^e}(\mathcal{P}_{\rm Hoch},\Lambda)}$. Hence
 these maps induce inverse linear isomorphisms at the cohomology level. Taking
 into account the identifications above we obtain that
 $$\begin{array}{rrcl}
   \ov{\omega_1}:&k(Q_1//B)&\longrightarrow&\Hom_k(\Lambda,\Lambda)\\
   &(a,\gamma)&\longmapsto&
    \begin{array}[t]{rcl}
     \Lambda&\rightarrow&\Lambda\\
     \varepsilon&\mapsto&\varepsilon^{(a,\gamma)}
    \end{array}\qquad\mbox{\rm and}\\
   \ov{\varsigma_1}:&\Hom_k(\Lambda,\Lambda)&\longrightarrow&k(Q_1//B)\\
   &\begin{array}[t]{rcl}
    f:\Lambda&\rightarrow&\Lambda\\
     \varepsilon&\mapsto&\sum_{\gamma\in B}\lambda_{\varepsilon,\gamma}\gamma
    \end{array}
   &\longmapsto&\sum_{a\in Q_1}\sum_{(a,\gamma)\in Q_1//B}
    \lambda_{a,\gamma}(a,\gamma)
 \end{array}$$
 induce inverse linear isomorphisms between ${\rm H}^1(\Lambda,\Lambda)=
 \Ker d_1/\im d_0$ and
 ${\rm H}^1(\Lambda,\Lambda)=\Ker\psi_1/\im\psi_0$.
 This allows us to transfer the Lie algebra structure of $\Ker d_1/\im d_0$ to
 $\Ker\psi_1/\im\psi_0$. Define on $k(Q_1//B)$ the bracket
 $$[(a,\gamma),(b,\varepsilon)]:=\ov{\varsigma_1}([\ov{\omega_1}(a,\gamma),
   \ov{\omega_1}(b,\varepsilon)])=(b,\varepsilon^{(a,\gamma)})-
   (a,\gamma^{(b,\varepsilon)})$$
 for all $(a,\gamma),(b,\varepsilon)\in Q_1//B$. It is easy to check that this
 is a Lie bracket. As $\ov{\omega}:\mathcal{C}_{\rm Hoch}\rightarrow
 \mathcal{C}_{\rm min}$ and $\ov{\varsigma}:\mathcal{C}_{\rm min}\rightarrow
 \mathcal{C}_{\rm Hoch}$ are maps of complexes, we conclude from the fact that 
 $\Ker d_1=\Der_k(\Lambda)$ is a Lie subalgebra of $\Hom_k(\Lambda,\Lambda)$, 
 that $\Ker\psi_1$ is a Lie subalgebra of $k(Q_1//B)$. In the same way we 
 deduce from the fact that $\im\psi_0=\Ad_k(\Lambda)$ is a Lie ideal of 
 $\Ker d_1$, that $\im\psi_0$ is a Lie ideal of $\Ker\psi_1$. By construction 
 the quotient Lie algebra ${\rm H}^1(\Lambda,\Lambda)=\Ker\psi_1/\im\psi_0$ is
 isomorphic to the Lie algebra ${\rm H}^1(\Lambda,\Lambda)=\Ker d_1/\im d_0$.
\end{Bew}

It is remarkable that Bardzell's minimal projective resolution and the 
resolution considered by Gerstenhaber and Schack in \cite{GS2} using 
$E$-relative Hochschild cohomology yield the same $1$-coboundaries and 
$1$-cocycles as the following Proposition and Corollary show:

\begin{Prop}
 The Lie algebras $\Ker\psi_1$ (with the bracket described in the 
 preceding Theorem) and $\Der_{E^e}(\Lambda)$ (with the canonical bracket) 
 are isomorphic.
\end{Prop}

\begin{Bew}
 Let $\sigma:\Lambda\rightarrow\Lambda$ be an element of 
 $\Der_{E^e}(\Lambda)=\Der_k(\Lambda)\cap\Hom_{E^e}(\Lambda,\Lambda)$.
 Since $\sigma$ is an $E$-bimodule morphism, the image of an arrow $a$ is a 
 linear combination of paths of $B$ parallel to $a$. Suppose
 $\sigma(a)=\sum_{(a,\gamma)\in Q_1//B}\lambda_{(a,\gamma)}\gamma$ for 
 $a\in Q_1$ with $\lambda_{(a,\gamma)}\in k$. The fact that $\sigma$ is a 
 well-defined derivation implies that we have for every path  
 $p=p_l\ldots p_2p_1$ in $Z$ of length $l$
 $$0=\sigma(p)=\sigma(p_l\ldots p_1)=
   \sum_{i=1}^lp_l\ldots p_{i+1}\sigma(p_i)p_{i-1}\ldots p_1=
   \sum_{(a,\gamma)\in Q_1//B}\lambda_{(a,\gamma)}p^{(a,\gamma)}$$
 Thus we have in $k(Z//B)$ the equality 
 $0=\sum_{p\in Z}\sum_{(a,\gamma)\in Q_1//B}
  \lambda_{(a,\gamma)}(p,p^{(a,\gamma)})$
 which is the same as to say that 
 $\sum_{(a,\gamma)\in Q_1//B}\lambda_{(a,\gamma)}(a,\gamma)$ is an element in
 $\Ker\psi_1$. This enables us to define the linear function
 $$\xi:\Der_{E^e}(\Lambda)\rightarrow\Ker\psi_1,\qquad
   \sigma\mapsto\sum_{(a,\gamma)\in Q_1//B}\lambda_{(a,\gamma)}(a,\gamma)$$
 On the other hand we can associate to every 
 $x=\sum_{(a,\gamma)\in Q_1//B}\lambda_{(a,\gamma)}(a,\gamma)$ in $\Ker\psi_1$
 an $E$-bimodule morphism $\sigma_x:\Lambda\rightarrow\Lambda$ by setting
 $\sigma_x(e):=0$ for every vertex $e\in Q_0$ and 
 $\sigma_x(c):=\sum_{(c,\gamma)\in Q_1//B}\lambda_{(c,\gamma)}\gamma$ for every
 arrow $c\in Q_1$. To make $\sigma_x$ a derivation we extend $\sigma_x$ to 
 paths $\varepsilon=b_m\ldots b_1$ of length $m\ge2$ by the formula
 $$\sigma_x(\varepsilon)=\sigma_x(b_m\ldots b_1)=
   \sum_{j=1}^mb_m\ldots b_{j+1}\sigma_x(b_j)b_{j-1}\ldots b_1=
   \sum_{(a,\gamma)\in Q_1//B}\lambda_{(a,\gamma)}\varepsilon^{(a,\gamma)}$$
 Thanks to $0=\psi_1(x)=\sum_{p\in Z}\sum_{(a,\gamma)\in Q_1//B}
 \lambda_{(a,\gamma)}(p,p^{(a,\gamma)})$ we have for every $p\in Z$ the 
 relation $0=\sum_{(a,\gamma)\in Q_1//B}\lambda_{(a,\gamma)}p^{(a,\gamma)}
 =\sigma_x(p)$ which shows that $\sigma_x$ is well-defined. This allows us to 
 define the linear function 
 $$\zeta:\Ker\psi_1\rightarrow\Der_{E^e}(\Lambda),\qquad x\mapsto\sigma_x$$
 It is clear that $\xi$ and $\zeta$ are inverse to each other which proves that 
 $\Der_{E^e}(\Lambda)$ and $\Ker\psi_1$ are isomorphic $k$-vector spaces.
 In order to check that $\zeta$ is a Lie algebra morphism let us fix
 $x=\sum_{(a,\gamma)\in Q_1//B}\lambda_{(a,\gamma)}(a,\gamma)$ and
 $y=\sum_{(b,\varepsilon)\in Q_1//B}\mu_{(b,\varepsilon)}(b,\varepsilon)$ in 
 $\Ker\psi_1$. For every arrow $c\in Q_1$ we have
 \begin{eqnarray*}
  [\sigma_x,\sigma_y](c)&=&(\sigma_x\circ\sigma_y-\sigma_y\circ\sigma_x)(c)\\
  &=&\sum_{(c,\varepsilon)\in Q_1//B}\mu_{(c,\varepsilon)}\sigma_x(\varepsilon)
  -\sum_{(c,\gamma)\in Q_1//B}\lambda_{(c,\gamma)}\sigma_y(\gamma)\\
  &=&\sum_{(c,\varepsilon),(a,\gamma)\in Q_1//B}\mu_{(c,\varepsilon)}
     \lambda_{(a,\gamma)}\varepsilon^{(a,\gamma)}-
  \sum_{(c,\gamma),(b,\varepsilon)\in Q_1//B}\lambda_{(c,\gamma)}
     \mu_{(b,\varepsilon)}\gamma^{(b,\varepsilon)}\\
  &=&\sigma_z(c)
 \end{eqnarray*}
 where $z=\sum_{(a,\gamma),(b,\varepsilon)\in Q_1//B}
 \lambda_{(a,\gamma)}\mu_{(b,\varepsilon)}
 ((b,\varepsilon^{(a,\gamma)})-(a,\gamma^{(b,\varepsilon)}))$. Theorem
 \ref{bracket} shows that $z$ is equal to the element $[x,y]$ of $\Ker\psi_1$.
 Therefore $\zeta$ is a Lie algebra morphism and $\Der_{E^e}(\Lambda)$ and 
 $\Ker\psi_1$ are isomorphic Lie algebras.
\end{Bew}

\begin{Kor}
\label{out1}
 The Lie ideal $\im\psi_0$ of $\Ker\psi_1$ and the Lie ideal 
 $\Ad_{E^e}(\Lambda)$ of $\Der_{E^e}(\Lambda)$ are isomorphic.
\end{Kor}

\begin{Bew}
 Let $y=\sum_{\gamma\in B}\lambda_\gamma\gamma\in\Lambda$ be such that $\ad(y)$
 is an $E$-bimodule morphism. Then we have for every vertex $e\in Q_0$
 $$0=eye-eye=e\,\ad(y)(e)\,e=\ad(y)(e)=ye-ey=\sum_{\gamma\in B}\lambda_\gamma
   (\delta_{s(\gamma),e}-\delta_{t(\gamma),e})\gamma$$
 This shows that every path $\gamma\in B$ satisfying $\lambda_\gamma\neq0$ is a
 cycle and so 
 $y=\sum_{(e,\gamma)\in Q_0//B}\lambda_\gamma\gamma$. For every 
 arrow $a\in Q_1$ we have
 $$\ad(y)(a)=ya-ay=\sum_{(e,\gamma)\in Q_0//B}\lambda_\gamma\gamma a-
   \sum_{(e,\gamma)\in Q_0//B}\lambda_\gamma a\gamma=
   \sum_{(e,\gamma)\in Q_0//B}\lambda_\gamma(\gamma a-a\gamma)$$
 We deduce from this that the image of $\ad(y)\in\Ad_{E^e}(\Lambda)$ by $\xi$ 
 is 
 $$\xi(\ad(y))=\sum_{(e,\gamma)\in Q_0//B}\lambda_\gamma
   (\sum_{\stackrel{a\in eQ_1}{\gamma a\in B}}(a,\gamma a)-
    \sum_{\stackrel{a\in Q_1e}{a\gamma\in B}}(a,a\gamma))=
   -\sum_{(e,\gamma)\in Q_0//B}\lambda_\gamma\psi_0(e,\gamma)$$
 and thus we obtain that the image of the Lie algebra $\Ad_{E^e}(\Lambda)$ by 
 the isomorphism of Lie algebras $\xi$ is $\im\psi_0$.
\end{Bew}
 
\begin{Bem}
 In general, the minimal projective resolution and the resolution considered by
 Gerstenhaber and Schack do not yield the same $n$-coboundaries and 
 $n$-cocycles: for example in the case of hereditary algebras, \ie 
 $Z=\emptyset$, we have for $n=2$ that $\Ker\psi_2=0$, because $k(Z//B)=0$, but
 the multiplication $m\in\Hom_{E^e}(\Lambda\otimes_E\Lambda,\Lambda)$ is a 
 nontrivial $2$-coboundary as well as a $2$-cocycle of the complex yielded by 
 Gerstenhaber and Schack's resolution. 
\end{Bem}

\section{The algebraic group of outer automorphisms of a monomial algebra}

In this section, we assume that the characteristic of the field $k$ is $0$. Our
aim is to relate  the description of the Lie algebra 
${\rm H}^1(\Lambda,\Lambda)=\Ker\psi_1/\im\psi_0$ obtained in the preceding
section to the algebraic groups which appear in Guil-Asensio and Saor\'{\i}n's 
study of the outer automorphisms. Denote by $\Aut(A)$ the algebraic 
group of all $k$-algebra automorphisms of a finite dimensional $k$-algebra $A$.
The group of inner automorphisms $\Inn(A)$ of $A$ is the image of the morphism
of algebraic groups $\varphi:A^\ast\rightarrow\Aut(A)$ given by 
$a\mapsto \iota_a=$ conjugation by $a$. Thus $\Inn(A)$ is a closed normal and
connected subgroup of $\Aut(A)$. The algebraic group of outer automorphisms
$\Out(A)$ is defined as the quotient $\Aut(A)/\Inn(A)$. Its identity component
$\Out(A)^\circ$ is $\Aut(A)^\circ/\Inn(A)$. Note that if $A$ is a basic algebra,
for instance the monomial algebra $\Lambda$, then the group $\Out(A)$ is
isomorphic to the Picard group $\Pic_k(A)$, that is the group of (isomorphism
types of) Morita autoequivalences of the category of left $A$-modules (see
\cite{Po}, p. 1860).

\begin{Prop}
\label{Out}
 Let $k$ be a field of characteristic $0$ and $A$ a finite dimensional
 $k$-algebra. The derivations $\Der_k(A)$ form the Lie algebra of the algebraic
 group $\Aut(A)$ and the inner derivations $\Ad_k(A)$ form the Lie algebra of
 the algebraic group $\Inn(A)$. The Lie algebra ${\rm H}^1(A,A)$ can be
 regarded as the Lie algebra of the algebraic group $\Out(A)$ or as the Lie
 algebra of its identity component $\Out(A)^\circ$.
\end{Prop}

\begin{Bew}
 In \cite{Hu} 13.2 it is shown that $\Der_k(A)$ is the Lie algebra of $\Aut(A)$.
 The differential ${\rm d}\varphi:A\rightarrow\Der_k(A)$ of the above-mentioned
 morphism $\varphi$ is given by $a\mapsto \ad(a)=$ inner derivation of $a$. 
 Since
 the field $k$ has characteristic $0$, the Lie algebra of the image of a morphism
 of algebraic groups is the image of the differential of this morphism. The
 construction of the quotient $\Aut(A)/\Inn(A)$ implies that the Lie algebra of
 the algebraic group $\Aut(A)/\Inn(A)=\Out(A)$ is 
 $\mathcal{L}(\Aut(A))/\mathcal{L}(\Inn(A))=\Der_k(A)/\Ad_k(A)={\rm H}^1(A,A)$
 (see \cite {Hu} 11.5 and 12).
\end{Bew}

The automorphism group of any finite dimensional algebra was studied by Pollack
in \cite{Po}. Guil-Asensio and Saor\'{\i}n worked on the outer automorphisms of
any finite dimensional algebra (see \cite{GAS1}). The case of the finite
dimensional monomial algebras was treated by them in their paper $\cite{GAS2}$.
We are going to follow their notations.
Define the algebraic group
$$H_\Lambda:=\{\sigma\in\Aut(\Lambda)\mid \sigma(e)=e\ \forall e\in Q_0\}$$
and denote by $H_\Lambda^\circ$ its identity component. According to 
\cite{Po} and to
Proposition 1.1 of \cite{GAS2} we have
$$\Out(\Lambda)^\circ\simeq H_\Lambda^\circ/H_\Lambda\cap\Inn(\Lambda)$$

\begin{Prop}
 Let $k$ be a field of characteristic $0$ and $\Lambda$ a finite dimensional
 monomial algebra. The Lie algebra of the algebraic group $H_\Lambda$ is the Lie
 algebra $\Der_{E^e}(\Lambda)\simeq\Ker\psi_1$ and
 $\Ad_{E^e}(\Lambda)\simeq\im\psi_0$ is the Lie algebra of
 $H_\Lambda\cap\Inn(\Lambda)$.
\end{Prop}

\begin{Bew}
 For every vertex 
 $e\in Q_0$ we write $G_e:=\{\sigma\in{\bf GL}(\Lambda)\mid\sigma(e)=e\}$ and 
 $\mathcal{G}_e:=\{\sigma\in\mathfrak{gl}(\Lambda)\mid\sigma(e)=0\}$. Thus we get the
 equality $H_\Lambda=\Aut(\Lambda)\cap\bigcap_{e\in Q_0}G_e$. Since we have
 assumed $\car k=0$, paragraph 13.2 of \cite{Hu} shows that
 $$\begin{array}{rcl}
   \mathcal{L}(H_\Lambda)&=&\mathcal{L}(\Aut(\Lambda))\cap\bigcap_{e\in Q_0}
   \mathcal{L}(G_e)
   \ =\ \Der_k(\Lambda)\cap\bigcap_{e\in Q_0}\mathcal{G}_e\\
   &=&\{\sigma\in\Der_k(\Lambda)\mid\sigma(e)=0\ \forall e\in Q_0\}
   \ =\ \Der_{E^e}(\Lambda)\qquad\quad{\rm and}\\
   &&\\
 \mathcal{L}(H_\Lambda\cap\Inn(\Lambda))
 &=&\mathcal{L}(H_\Lambda)\cap
   \mathcal{L}(\Inn(\Lambda))\ =\ \Der_{E^e}(\Lambda)\cap\Ad_k(\Lambda)
   \ =\ \Ad_{E^e}(\Lambda)
 \end{array}$$
\end{Bew}

 The closed unipotent and connected subgroup
$$\hat{E}:=\varphi(E^\ast)=\{\iota_a\in\Inn(\Lambda)\mid a\in E^\ast\}$$ 
of $H_\Lambda$ is isomorphic
as a group to the group of acyclic characters ${\rm Ch}(Q,k)$ appearing in 
Guil-Asensio and Saor\'{\i}n's paper (see definition 7 in \cite {GAS1}). The 
group
$$\Inn^\ast(\Lambda):=\varphi(1+r)=\{\iota_a\in\Inn(\Lambda)\mid
 \exists x\in r\ {\rm such\ that}\ a=1+x\}$$
is a closed unipotent and connected subgroup of $\Inn(\Lambda)$. We have
$H_\Lambda\cap\Inn(\Lambda)=(H_\Lambda\cap\Inn^\ast(\Lambda))\times\hat{E}$.
Define the algebraic groups 
$$U_\Lambda:=
  \{\sigma\in\Aut(\Lambda)\mid\sigma(a)\equiv a\mod r^2\ \forall a\in Q_1\}
  \qquad {\rm and}$$
$$V_\Lambda:=\{\sigma\in H_\Lambda\mid\sigma(a)=\sum_{(a,b)\in Q_1//Q_1}
 \lambda_{(a,b)}b\ \ \forall a\in Q_1\}$$
According to Proposition 1.1 of \cite{GAS2} we have 
$H_\Lambda^\circ=(H_\Lambda\cap U_\Lambda)\rtimes V_\Lambda^\circ$ and therefore
$\Out(\Lambda)^\circ$ can be described as follows:

\begin{Prop}
\label{semidirect}
 Let $k$ be a field of characteristic $0$ and $\Lambda$ a finite dimensional
 monomial $k$-algebra. Then the identity component $\Out(\Lambda)^\circ$ of 
 the algebraic group of
 outer automorphisms is the semidirect product
 $$\Out(\Lambda)^\circ=
   \frac{H_\Lambda\cap U_\Lambda}{H_\Lambda\cap\Inn^\ast(\Lambda)}\rtimes
   \frac{V_\Lambda^\circ}{\hat{E}}$$
\end{Prop}

Since there is an inclusion preserving $1-1$ correspondence between the
collection of closed connected subgroups of $\Out(\Lambda)^\circ$ (resp.
$H_\Lambda^\circ$) and the collection of their Lie algebras, regarded as
subalgebras of ${\rm H}^1(\Lambda,\Lambda)=\Ker\psi_1/\im\psi_0$ 
(resp. $\Ker\psi_1$) we are interested in identifying the subalgebras of 
${\rm H}^1(\Lambda,\Lambda)=\Ker\psi_1/\im\psi_0$ (resp. $\Ker\psi_1$) 
corresponding to the algebraic groups 
$\frac{H_\Lambda\cap U_\Lambda}{H_\Lambda\cap\Inn^\ast(\Lambda)}$ and 
$\frac{V_\Lambda^\circ}{\hat{E}}$ (resp. $H_\Lambda\cap U_\Lambda$,
$H_\Lambda\cap\Inn^\ast(\Lambda)$, $V_\Lambda^\circ$ and $\hat{E}$).
We have the following dictionary:

\begin{Prop}
 Let $k$ be a field of characteristic $0$ and $\Lambda$ a finite dimensional
 monomial $k$-algebra.
 \begin{enumerate}
  \item The Lie algebra of the closed connected normal subgroup 
        $H_\Lambda\cap U_\Lambda$ of $H_\Lambda^\circ$ is the Lie ideal 
        $k(Q_1//B-(Q_0\cup Q_1))\cap\Ker\psi_1$ of the Lie algebra $\Ker\psi_1$.
  \item The Lie algebra of the closed connected  subgroup 
        $V_\Lambda^\circ$ of $H_\Lambda^\circ$ is the Lie subalgebra 
        $k(Q_1//Q_1)\cap\Ker\psi_1$ of the Lie algebra $\Ker\psi_1$.
  \item The Lie algebra of the closed connected normal subgroup 
        $H_\Lambda\cap\Inn^\ast(\Lambda)$ of \linebreak
    $H_\Lambda\cap\Inn(\Lambda)$ is 
    the Lie ideal of $\im\psi_0$ generated by the elements
    $\sum_{a\in Q_1e}(a,a\gamma)-\sum_{a\in eQ_1}(a,\gamma a)$ where
    $(e,\gamma)\in Q_0//(B-Q_0)$.
  \item The Lie algebra of the closed connected subgroup 
        $\hat{E}$ of $H_\Lambda\cap\Inn(\Lambda)$ is 
    the Lie subalgebra of $\im\psi_0$ generated by the elements
    $\sum_{a\in Q_1e}(a,a)-\sum_{a\in eQ_1}(a,a)$ where $e\in Q_0$.
 \end{enumerate}
\end{Prop}

\begin{Bew}
 The derivation of the morphism of algebraic groups 
 $$\begin{array}{rrcl}
 \varepsilon_\Lambda:&\Aut(\Lambda)&\longrightarrow&
 {\bf GL}(r/r^2)\simeq{\bf GL}(kQ_1)\\
 &\sigma&\longmapsto&
   r/r^2\rightarrow r/r^2,\ \ 
   x\mod r^2\mapsto\sigma(x)\mod r^2\\
 \end{array}$$
 is given by
 $$\begin{array}{rrcl}
 {\rm d}\varepsilon_\Lambda:&\Der_k(\Lambda)&\longrightarrow&
 \mathfrak{gl}(r/r^2)\simeq\mathfrak{gl}(kQ_1)\\
 &\sigma&\longmapsto&
   r/r^2\rightarrow r/r^2,\
   x\mod r^2\mapsto\sigma(x)\mod r^2\\
 \end{array}$$
 Recall that we have shown in Proposition 1.8 that
 $$\begin{array}{rccl}
   \xi:&\Der_{E^e}(\Lambda)&\longrightarrow&\Ker\psi_1\\
   &\begin{array}[t]{rrcl}
    \sigma:&\Lambda&\rightarrow&\Lambda\\
    &a\in Q_1&\mapsto&\sum_{(a,\gamma)\in Q_1//B}\lambda_{(a,\gamma)}(a,\gamma)
    \end{array}
   &\longmapsto&\sum_{(a,\gamma)\in Q_1//B}\lambda_{(a,\gamma)}(a,\gamma)
   \end{array}$$
  is a Lie algebra isomorphism.
  
 {\it (i)}\ : Since $U_\Lambda$ is the kernel of $\varepsilon_\Lambda$ we 
 obtain, thanks to the assumption $\car k=0$,
 $$\mathcal{L}(H_\Lambda\cap U_\Lambda)=
   \mathcal{L}(H_\Lambda)\cap\Ker{\rm d}\varepsilon_\Lambda=
   \{\sigma\in\Der_{E^e}(\Lambda)\mid\sigma(a)\in r^2\ \forall a\in Q_1\}$$
 The fact $\xi(\mathcal{L}(H_\Lambda\cap U_\Lambda))=k(Q_1//B-(Q_0\cup Q_1))\cap
 \Ker\psi_1$ finishes the proof.
 
 {\it (ii)}\ : From $V_\Lambda\simeq\im\varepsilon_\Lambda\mid_{H_\Lambda}$
 (see Lemma 22 in \cite{GAS1}) and $\car k=0$ we deduce
 $$\begin{array}{rcl}
   \mathcal{L}(V_\Lambda)&=&\mathcal{L}(V_\Lambda^\circ)\ =\ 
 \im{\rm d}\varepsilon_\Lambda\mid_{\mathcal{L}(H_\Lambda)}\\
 &=&\{\sigma\in\Der_{E^e}(\Lambda)\mid
 \sigma(a)=\sum_{(a,b)\in Q_1//Q_1}\lambda_{(a,b)}(a,b)\ \forall a\in Q_1\}
 \end{array}$$
 The equality $\xi(\mathcal{L}(V_\Lambda^\circ))=k(Q_1//Q_1)\cap\Ker\psi_1$
 shows that we are done.
 
 {\it (iii)}\ : We have
 $$\mathcal{L}(H_\Lambda\cap\Inn^\ast(\Lambda))= 
   \{\ad(1+x)\in\Ad_{E^e}(\Lambda)\mid x\in r\}=
   \{\ad(x)\in\Ad_{E^e}(\Lambda)\mid x\in r\}$$
 Let $x$ be an element of the radical $r$ such that $\ad(x)$ is an $E$-bimodule
 morphism. The proof of Corollary \ref{out1} shows that $x$ is a linear
 combination of oriented cycles of length $\ge1$, so
 $x=\sum_{(e,\gamma)\in Q_0//B-Q_0}\lambda_\gamma\gamma$ with 
 $\lambda_\gamma\in k$. Since $\xi(\ad(x))=\sum_{(e,\gamma)\in Q_0//B-Q_0}
 \lambda_\gamma(\sum_{\stackrel{a\in eQ_1}{\gamma a\in B}}(a,\gamma a)-
 \sum_{\stackrel{a\in Q_1e}{a\gamma\in B}}(a,a\gamma))$, we are done.
 
 {\it (iv)}\ : It is obvious that the Lie algebra of $\hat{E}$ is generated by
 the inner derivations $\ad(e)$ where $e\in Q_0$. Since 
 $\xi(\ad(e))=\sum_{a\in Q_1e}(a,a)-\sum_{a\in eQ_1}(a,a)$ the proof is
 finished.
\end{Bew}

Using the notations which will be introduced at the beginning of the next
section we get:

\begin{Kor}
\label{dico}
 Let $k$ be a field of characteristic $0$ and $\Lambda$ a finite dimensional
 monomial $k$-algebra. The Lie algebra of the closed connected normal subgroup 
 $\frac{H_\Lambda\cap U_\Lambda}{H_\Lambda\cap\Inn^\ast(\Lambda)}$ of the
 connected algebraic group $\Out(\Lambda)^\circ$ is the Lie ideal
 $\oplus_{i\ge1}L_i$ of ${\rm H}^1(\Lambda,\Lambda)$ and the Lie algebra of the
 closed connected subgroup $\frac{V_\Lambda^\circ}{\hat{E}}$ of
 $\Out(\Lambda)^\circ$ is the Lie subalgebra $L_0$ of 
 ${\rm H}^1(\Lambda,\Lambda)$.
\end{Kor}

\section{The Lie algebra ${\rm H}^1(\Lambda,\Lambda)$ of a monomial algebra}

Since Hochschild cohomology is additive and since its Lie algebra structure 
follows this additive decomposition we will assume henceforth that the quiver 
$Q$ is connected. For the study of the Lie algebra ${\rm H}^1(\Lambda,\Lambda)$
of the monomial algebra $\Lambda=kQ/\langle Z\rangle$ we will use the description which we obtained in Theorem \ref{bracket}. For every element 
$x\in\Ker\psi_1$ we will also denote its class in 
${\rm H}^1(\Lambda,\Lambda)=\Ker\psi_1/\im\psi_0$ by $x$. 

If $(a,\gamma)\in Q_1//B_n$ and $(b,\varepsilon)\in Q_1//B_m$, the formula we 
have obtained shows that $[(a,\gamma),(b,\varepsilon)]$ is an element of
$k(Q_1//B_{n+m-1})$. Thus, we have a graduation on the Lie algebra 
$k(Q_1//B)=\oplus_{i\in\N}k(Q_1//B_i)$ by considering that the elements of 
$k(Q_1//B_i)$ have degree $i-1$ for all $i\in\N$. It is clear that the Lie 
subalgebra $\Ker\psi_1$ of $k(Q_1//B)$ preserves this graduation and that 
$\im\psi_0$ is a graded ideal. Therefore, the Lie algebra 
${\rm H}^1(\Lambda,\Lambda)=\Ker\psi_1/\im\psi_0$ has also a graduation. If we
set
$$L_{-1}:=k(Q_1//Q_0)\cap\Ker\psi_1$$
$$L_0:=k(Q_1//Q_1)\cap\Ker\psi_1/\langle\sum_{a\in Q_1e}(a,a)-
  \sum_{a\in eQ_1}(a,a)\mid e\in Q_0\rangle\qquad{\rm and}$$
$$L_i:=k(Q_1//B_{i+1})\cap\Ker\psi_1/\langle
  \sum_{\stackrel{a\in Q_1e}{\gamma a\in B}}(a,\gamma a)-
  \sum_{\stackrel{a\in eQ_1}{a\gamma\in B}}(a,a\gamma)\mid
  (e,\gamma)\in Q_0//Q_i\rangle$$
for all $i\ge1$, $i\in\N$, we obtain 
${\rm H}^1(\Lambda,\Lambda)=\oplus_{i\ge-1}L_i$ and $[L_i,L_j]\subset L_{i+j}$
 for all $i,j\ge-1$ where $L_{-2}=0$. 

\begin{Lem}
 $L_{-1}$ equals $0$ if and only if there exists for every loop 
 $(a,e)\in Q_1//Q_0$ a path $p$ in $Z$ such that
 $p^{(a,e)}\neq0$.
\end{Lem} 

\begin{Bew}
 Clear, since $\psi_1(a,e)=\sum_{p\in Z}(p,p^{(a,e)})$
 for every loop $(a,e)\in Q_1//Q_0$.
\end{Bew}

We study first the subalgebra $L_{-1}$ of ${\rm H}^1(\Lambda,\Lambda)$.

\begin{Prop}
\label{L-1}
 Each of the following conditions implies $L_{-1}=0$:
 \begin{enumerate}
  \item The quiver $Q$ does not have a loop.
  \item For every loop $(a,e)\in Q_1//Q_0$ of $Q$ the characteristic of the 
    field $k$ does not divide the integer $m\ge2$ for which $a^m\in Z$ and 
    $a^{m-1}\in B$.
  \item The characteristic of the field $k$ is equal to $0$. 
  \item $\Lambda$ is a truncated quiver algebra $kQ/\langle Q_m\rangle$
    for a quiver $Q$ different from the loop and for $m\ge2$.
  \item $Q$ is the loop quiver and $\Lambda=kQ/\langle Q_m\rangle$ is a 
    truncated quiver algebra such that the characteristic of $k$ does not divide
    $m\ge2$.
 \end{enumerate}
\end{Prop}

\begin{Bew}
 {\it (i)}\ : clear

 {\it (ii)}\ : Let $(a,e)\in Q_1//Q_0$ be a loop of $Q$. Since $\Lambda$ is 
 finite dimensional, there exists an integer $m\ge2$ such that $p:=a^m\in Z$ 
 and $a^{m-1}\in B$. If the characteristic of $k$ does not divide $m$, then
 $p^{(a,e)}=ma^{m-1}$ is different from $0$ and so 
 $ma^{m-1}\not\in\langle Z\rangle$ \ie $p^{(a,e)}\neq0$.

 {\it (iii)}\ : Clear, because {\it (iii)} implies {\it (ii)}.

 {\it (iv)}\ : Using {\it (ii)} we can suppose that $\car k$ divides $m$. Let
 $(a,e)\in Q_1//Q_0$ be a loop of the quiver $Q$. By assumption the connected 
 quiver itself is not a loop, therefore there exists an arrow $b\in Q_1$ 
 different from $a$ such that $ba^{m-1}$ or $a^{m-1}b\in Z=Q_m$. The fact that
 $\car k$ divides $m$ implies that the characteristic of the field $k$ does not
 divide $m-1$ and thus 
 $0\neq(m-1)ba^{m-2}\not\in\langle Z\rangle=
 \langle Q_m\rangle$ or $0\neq(m-1)a^{m-2}b\not\in\langle Z\rangle=
 \langle Q_m\rangle$.

 {\it (v)}\ : Clear, because {\it (v)} implies {\it (ii)}.
\end{Bew}

Before we study the case where $L_{-1}$ equals zero we consider the following
exceptional case:

\begin{Prop}
\label{loop}
 Let $Q$ be the loop having a vertex $e$ and an arrow $a$. Let 
 $Z=\{a^m\}$ be
 such that the characteristic $p$ of the field $k$ divides the integer $m\ge2$
 and let $\Lambda=kQ/\langle Z\rangle$ be the monomial algebra associated to 
 $Z$. The following conditions are equivalent:
 \begin{enumerate}
  \item The Lie algebra ${\rm H}^1(\Lambda,\Lambda)$ is simple.
  \item The Lie algebra ${\rm H}^1(\Lambda,\Lambda)$ is semisimple.
  \item The integer $m\ge2$ is equal to the characteristic $p$ of the field $k$
        and $p>2$.
  \item ${\rm H}^1(\Lambda,\Lambda)$ is isomorphic to the Witt Lie algebra             
  $W(1,1):=\Der(k[X]/(X^p))$ and $p>2$.
 \end{enumerate}
\end{Prop}

\begin{Bew}
 {\it (i)}\ $\Rightarrow$\ {\it (ii)}\ : trivial
 
 {\it (ii)}\ $\Rightarrow$\ {\it (iii)}\ : The assumption $\car k=p$ implies 
 $[L_{-1},L_{p-1}]=0$ and thus, 
 $L_{p-1}\oplus\cdots\oplus L_{m-2}$ is a solvable
 Lie ideal of the Lie algebra 
 ${\rm H}^1(\Lambda,\Lambda)=\bigoplus_{i=-1}^{m-2}L_i$. Since 
 ${\rm H}^1(\Lambda,\Lambda)$ is semisimple, we deduce $m=p$. We have $p>2$, 
 because there is no semisimple Lie algebra of dimension $2$.

 {\it (iii)}\ $\Rightarrow$\ {\it (iv)}\ : Since $\Lambda$ is isomorphic to the 
 commutative algebra $k[X]/(X^m)$, the Lie algebra ${\rm H}^1(\Lambda,\Lambda)$
 is isomorphic to the Lie algebra of derivations $\Der(k[X]/(X^m))$.

 {\it (iv)}\ $\Rightarrow$\ {\it (i)}\ : If $\car k>2$, then the Witt Lie algebra 
 $W(1,1):=\Der(k[X]/(X^m))$ is one of the non classical simple Lie algebras.
\end{Bew}


Proposition \ref{L-1} shows that the case where $L_{-1}$ is different from
zero is quite exceptional. We will assume henceforth that $L_{-1}=0$. In that 
case $\bigoplus_{i\ge1}L_i$ is a solvable Lie ideal of 
${\rm H}^1(\Lambda,\Lambda)$ since ${\rm H}^1(\Lambda,\Lambda)$ is finite 
dimensional. It is obvious that $L_0$ is a Lie subalgebra of 
${\rm H}^1(\Lambda,\Lambda)$ whose bracket is
$$[(a,c),(b,d)]=\delta_{a,d}(b,c)-\delta_{b,c}(a,d)$$
for all $(a,c),(b,d)\in L_0$. It follows that we have
$\Rad{\rm H}^1(\Lambda,\Lambda)=\Rad L_0\oplus\bigoplus_{i\ge1}L_i$ and
${\rm H}^1(\Lambda,\Lambda)/\Rad{\rm H}^1(\Lambda,\Lambda)=L_0/\Rad L_0$ 
where $\Rad{\rm H}^1(\Lambda,\Lambda)$ (resp. $\Rad L_0$) denotes the radical
of ${\rm H}^1(\Lambda,\Lambda)$ (resp. $L_0$). As a consequence the study of the Lie algebra
${\rm H}^1(\Lambda,\Lambda)$ can be often reduced to the study of the Lie 
algebra $L_0$. 

We recall a few definitions introduced by Guil-Asensio and Saor\'{\i}n 
(see 2.3 in \cite{GAS2} and 25 in \cite{GAS1}) for the convenience of the 
reader which will be useful in the following. 

\begin{Def}
 Let $(a,b)$ be a couple of parallel arrows. We shall say that the ideal 
 $\langle Z\rangle$ of the algebra $kQ$ is $(a,b)$-saturated, if for every 
 path $p$ of $Z$ we have $p^{(b,a)}=0$. This is denoted by 
 $a\le_{\langle Z\rangle}b$. The ideal $\langle Z\rangle$ is called completely
 saturated if it is $(a,b)$-saturated for all $(a,b)\in Q_1//Q_1$.
\end{Def}

\begin{Beme}
\ \begin{enumerate}
  \item For every class of parallel arrows $\ov{\alpha}$ the relation 
   $\le_{\langle Z\rangle}$ on $\ov{\alpha}$ is reflexive and transitive.
  \item For $p\in Z$ and $(a,b)\in Q_1//Q_1$ we have $p^{(b,a)}=0$ if and only 
   if each term of the sum is zero (see \ref{distinct}), \ie each replacement 
   of one appearence of $b$ in $p$ by $a$ is a path in $Z$.
 \end{enumerate}
\end{Beme}

The hereditary algebra where $Z=\emptyset$ and the truncated algebras where 
$Z=Q_m$, $m\ge2$, are examples of monomial algebras whose ideal 
$\langle Z\rangle$ is completely saturated. Locateli studied the Hochschild 
cohomology of truncated algebras in \cite{Lo} and Cibils treated the case of 
radical square zero algebras (see \cite{Ci}).

It is clear that the relation of parallelism is an equivalence relation on 
$Q_1$. Denote by $Q_1//$ its set of equivalence classes. We shall call a class
of parallel arrows non trivial if it contains at least two arrows. It is easy 
to check that we have on every class of parallel arrows 
$\ov{\alpha}=\{\alpha_1,\ldots,\alpha_n\}\in Q_1//$ the following equivalence 
relation
$$\alpha_i\approx\alpha_j\quad :=\quad \alpha_i\le_{\langle Z\rangle}\alpha_j
                          \ {\rm and}\ \alpha_j\le_{\langle Z\rangle}\alpha_i$$
We denote by $\ov{\alpha}/\approx$ its set of equivalence classes. For 
$R,S\in\ov{\alpha}/\approx$ we write $R\le_{\langle Z\rangle}S$ if there
exist arrows $a\in R$ and $b\in S$ such that $a\le_{\langle Z\rangle}b$ which 
is equivalent to saying that $a\le_{\langle Z\rangle}b$ for all arrows $a\in R$ 
and $b\in S$. Note that the relation $\le_{\langle Z\rangle}$ on 
$\ov{\alpha}/\approx$ is an order relation.

The notion of saturation can be reformulated in terms of the 
Lie algebra structure:

\begin{Bem}
\label{saturated}
 Let $(a,b)$ be a couple of parallel arrows. The ideal $\langle Z\rangle$ is 
 $(a,b)$-saturated if and only if $(b,a)$ is an element of the Lie algebra 
 $L_0$. Thus the parallel arrows $a$ and $b$ are equivalent if and only if
 $(a,b)\in L_0$ and $(b,a)\in L_0$. The ideal $\langle Z\rangle$ is completely 
 saturated if and only if every two parallel arrows are equivalent.
\end{Bem}

\begin{Prop}
\label{basis}
 \ \\A basis $\mathcal{B}$ of the Lie algebra $L_0$ is given by the union of 
 the following sets:
 \begin{enumerate}
  \item all the couples $(a,b)\in L_0$ such that the parallel arrows $a$ and 
        $b$ are different
  \item for every class of parallel arrows 
        $\ov{\alpha}=\{\alpha_1,\ldots,\alpha_n\}\in Q_1//$, all the elements
        $(\alpha_i,\alpha_i)\in L_0$ such that $i<n$
  \item $\vert Q_1//\vert-\vert Q_0\vert+1$ linearly independent elements $(c,c)\in L_0$ 
        different from those in {\it (ii)}.
 \end{enumerate}
\end{Prop}

\begin{Bew}
 An element $\sum_{(a,b)\in Q_1//Q_1}\lambda_{(a,b)}(a,b)\in k(Q_1//Q_1)$, with
 $\lambda_{(a,b)}\in k$, is contained in $\Ker\psi_1$ if and only if
 $$0=\sum_{(a,b)\in Q_1//Q_1}\lambda_{(a,b)}\psi_1(a,b)=
   \sum_{\stackrel{(a,b)\in Q_1//Q_1}{a\neq b}}\lambda_{(a,b)}\sum_{p\in Z}(p,p^{(a,b)})$$
 First we recall that for a couple of distinct arrows $(a,b)$, the non zero
 terms of the sum defining $p^{(a,b)}$ are all distinct (see \ref{distinct}), 
 forming a subset of the basis $Z//B$. Second we check the following: if
 $a\neq b$ and if $(p,p^{(a,b)})=(q,q^{(c,d)})$ then $p=q$ and if 
 $p^{(a,b)}\neq 0$ then $(a,b)=(c,d)$. Indeed let $p=p_n\ldots p_2p_1$ be a
 path of $Z$ of length $l$ and consider the summand 
 $p_n\ldots p_{i+1}bp_{i-1}\ldots p_1\in B$
 of $p^{(a,b)}$ obtained by replacing the arrow $p_i=a$ by $b\neq a$. 
 The fact $(p,p^{(a,b)})=(p,p^{(c,d)})\neq0$ implies that there exists 
 $j\in\{1,\ldots,n\}$ such that
 $p_n\ldots p_{j+1}cp_{j-1}\ldots p_{i+1}bp_{i-1}\ldots p_1=
  p_n\ldots p_{j+1}dp_{j-1}\ldots p_{i+1}ap_{i-1}\ldots p_1$.
 Since $a\neq b$, we have $i=j$ as well as $a=d$ and $c=b$. This shows that if
 $a\neq b$ then the non zero elements $(p,p^{(a,b)})$ of $k(Z//B)$ are linearly
 independent. Therefore we can choose a basis of $k(Q_1//Q_1)\cap\Ker\psi_1$ 
 in $Q_1//Q_1$:
 take all elements $(a,a)\in Q_1//Q_1$ and all elements 
 $(a,b)\in Q_1//Q_1$, $a\neq b$, such that the ideal $\langle Z\rangle$ is 
 $(b,a)$-saturated. We have for every vertex $e\in Q_0$
 in $L_0=k(Q_1//Q_1)\cap\Ker\psi_1/k(Q_1//Q_1)\cap\im\psi_0$ the relation 
 $\sum_{a\in Q_1e}(a,a)-\sum_{a\in eQ_1}(a,a)=0$. Since $\vert Q_0\vert-1$ of 
 those $\vert Q_0\vert$ relations are linearly
 independent, we see that $\dim k(Q_1//Q_1)\cap\im\psi_0=\vert Q_1\vert-\vert Q_0\vert+1$.
 Thus we get a basis of $L_0$ if we impose the conditions ${\it (ii)}$ and 
 ${\it (iii)}$.
\end{Bew}

Let $\mathcal{B}$ be a basis of $L_0$ as described in the preceding Proposition. 

\begin{Bem}
\label{full}
 If the quiver $Q$ has a (non oriented or oriented) cycle, then there exists at least
 one class of parallel arrows $\ov{\alpha}=\{\alpha_1,\ldots,\alpha_n\}$, $n\ge1$, such that 
 $\mathcal{B}$ contains all the elements $(\alpha_i,\alpha_i)\in L_0$ , $1\le i\le n$.
\end{Bem}

For every class of parallel arrows
$\ov{\alpha}\in Q_1//$ we denote by $L_0^{\ov{\alpha}}$ the Lie ideal of $L_0$ 
generated by the elements $(a,b)\in\mathcal{B}$ such that $a,b\in\ov{\alpha}$.
Obviously the Lie algebra $L_0$ is the product of these Lie algebras:
$$L_0=\prod_{\ov{\alpha}\in Q_1//}L_0^{\ov{\alpha}}$$
To study the radical of this Lie algebra we need the following Lemma:

 \begin{Lem}
  Let $J\neq0$ be a Lie ideal of $L_0^{\ov{\alpha}}$ generated by elements
  $(\alpha_i,\alpha_j)\in L_0^{\ov{\alpha}}$ such that
  $(\alpha_j,\alpha_i)\not\in L_0^{\ov{\alpha}}$. Then $\dim_k[J,J]<\dim_kJ$.
 \end{Lem}
 
 \begin{Bew}
  Let $(\alpha_i,\alpha_j)$ be any element of $J$. If 
  $(\alpha_i,\alpha_j)\not\in[J,J]$ then nothing is to show. So suppose
  $(\alpha_i,\alpha_j)\in[J,J]$. This is the case if and only if there 
  exists an arrow $\alpha_l\in\ov{\alpha}$, 
  $\alpha_i\neq\alpha_l\neq\alpha_j$, such that $(\alpha_i,\alpha_l)\in J$ and
  $(\alpha_l,\alpha_j)\in J$. If $(\alpha_i,\alpha_l)\not\in[J,J]$ then we are
  done. If not we start again. Suppose that we have 
  $(\alpha_i,\alpha_l)\in[J,J]$. Thus there exists an arrow 
  $\alpha_k\in\ov{\alpha}-\{\alpha_i,\alpha_l\}$ such that 
  $(\alpha_i,\alpha_k)\in J$ and $(\alpha_k,\alpha_l)\in J$. The fact 
  $(\alpha_l,\alpha_j)\in J$ implies $(\alpha_j,\alpha_l)\not\in J$ and so
  $\alpha_k\neq\alpha_j$. If $(\alpha_i,\alpha_k)\not\in[J,J]$ then we are 
  done. If not we start again. Since $\ov{\alpha}$ is a finite set, this 
  procedure stops after a finite number of steps.
 \end{Bew}


\begin{Thm}
\label{radical}
 The radical of the Lie algebra $L_0$ is generated as a $k$-vector space by 
  the following elements of $L_0$: for every class of parallel arrows 
  $\ov{\alpha}=\{\alpha_1,\ldots,\alpha_n\}$
  \begin{itemize}
   \item $\sum_{\alpha_i\in S}(\alpha_i,\alpha_i)$ for every equivalence class
         $S\in\ov{\alpha}/\approx$
   \item $(\alpha_i,\alpha_j)$, $i\neq j$, if $(\alpha_i,\alpha_j)\in L_0$ and
         $(\alpha_j,\alpha_i)\not\in L_0$
   \item $(\alpha_{i_1 },\alpha_{i_2})$, $(\alpha_{i_2},\alpha_{i_1})$ and
         $(\alpha_{i_1},\alpha_{i_1})$ for all 
         $S=\{\alpha_{i_1},\alpha_{i_2}\}\in\ov{\alpha}/\approx$ if $\car k=2$
  \end{itemize}
\end{Thm}

\begin{Bew}
 Let $I$ be the $k$-vector space generated by the above-described elements. If 
 we define $I_{\ov{\alpha}}:=I\cap L_0^{\ov{\alpha}}$ we obtain 
 $I=\prod_{\ov{\alpha}\in Q_1//}I_{\ov{\alpha}}$. 
 Let $\ov{\alpha}=\{\alpha_1,\ldots,\alpha_n\}$ be a class of parallel arrows 
 of $Q$. It is easy to check using the definition of the bracket on $L_0$, 
 that $I_{\ov{\alpha}}$ is a Lie ideal of $L_0^{\ov{\alpha}}$ and that
 $I_{\ov{\alpha}}^{(3)}:=[I_{\ov{\alpha}}^{(2)},I_{\ov{\alpha}}^{(2)}]$ only
 contains elements $(\alpha_i,\alpha_j)\in L_0^{\ov{\alpha}}$ such that 
 $(\alpha_i,\alpha_j)\not\in L_0^{\ov{\alpha}}$. From the fact that
 $\dim_kI_{\ov{\alpha}}<\infty$ we deduce successively using the preceding 
 Lemma that there exists an $l\in\N$ such that 
 $I_{\ov{\alpha}}^{(l)}:=[I_{\ov{\alpha}}^{(l-1)},I_{\ov{\alpha}}^{(l-1)}]=0$.
 Thus $I_{\ov{\alpha}}$ is a solvable Lie ideal of $L_0^{\ov{\alpha}}$ and 
 we are done if we show that $L_0^{\ov{\alpha}}/I_{\ov{\alpha}}$ is a 
 semisimple Lie algebra or $0$.   
 From the relations in $L_0^{\ov{\alpha}}/I_{\ov{\alpha}}$ follows that we obtain a 
 basis of $L_0^{\ov{\alpha}}/I_{\ov{\alpha}}$ by 
 taking all elements $(\alpha_i,\alpha_j)$, where $\alpha_i\neq\alpha_j$ 
 and  $\alpha_i,\alpha_j\in S$, and $\vert S\vert-1$ elements 
 $(\alpha_i,\alpha_i)$ where $\alpha_i\in S$ for every 
 $S\in\ov{\alpha}/\approx$, with $\vert S\vert>2$ in case $\car k=2$. For every 
 $S\in\ov{\alpha}/\approx$ we denote by ${^SL}_0^{\ov{\alpha}}/I_{\ov{\alpha}}$
 the Lie ideal generated by the elements 
 $(\alpha_i,\alpha_j)\in L_0^{\ov{\alpha}}$ such that $\alpha_i,\alpha_j\in S$.
 Then we have $L_0^{\ov{\alpha}}/I_{\ov{\alpha}}=\prod_{S\in\ov{\alpha}/\approx}
 {^SL}_0^{\ov{\alpha}}/I_{\ov{\alpha}}$.
 Let $S=\{\alpha_{i_1},\ldots,\alpha_{i_m}\}$ be an equivalence class in 
 $\ov{\alpha}/\approx$, with $m\ge2$ if $\car k\neq2$ and $m\neq2$ if $\car k=2$. The elements 
 $(\alpha_{i_p},\alpha_{i_q})$, $i_p\neq i_q$, and 
 $(\alpha_{i_p},\alpha_{i_p})$, $i_p\neq i_m$, form a basis of 
 ${^SL}_0^{\ov{\alpha}}/I_{\ov{\alpha}}$ and we have the relation
 $\sum_{\alpha_{i_p}\in S}(\alpha_{i_p},\alpha_{i_p})=0$. As the bracket on
 ${^SL}_0^{\ov{\alpha}}/I_{\ov{\alpha}}$ is given by 
 $[(\alpha_{i_p},\alpha_{i_q}),(\alpha_{i_k},\alpha_{i_l})]=
  \delta_{i_q,i_k}(\alpha_{i_p},\alpha_{i_l})-
  \delta_{i_p,i_l}(\alpha_{i_k},\alpha_{i_q})$ we see that 
 ${^SL}_0^{\ov{\alpha}}/I_{\ov{\alpha}}\rightarrow\mathfrak{pgl}(m,k):=\mathfrak{gl}(m,k)/k1$ given by
 $(\alpha_{i_p},\alpha_{i_q})\mapsto e_{pq}\mod k1$ is a Lie algebra 
 isomorphism. Since $\mathfrak{pgl}(m,k)$ is a semisimple Lie algebra 
 if $m\ge2$ when $\car k\neq2$ and if $m>2$ when $\car k=2$, we conclude that 
 $L_0^{\ov{\alpha}}/I_{\ov{\alpha}}$ is a semisimple Lie algebra.
\end{Bew}

\begin{Kor}
\label{pgl}
 The semisimple Lie algebra 
 ${\rm H}^1(\Lambda,\Lambda)/{\rm Rad}\,{\rm H}^1(\Lambda,\Lambda)=
 L_0/{\rm Rad}\, L_0$ is
 the product of Lie algebras having a factor $\mathfrak{pgl}(\vert S\vert,k):=
 \mathfrak{gl}(\vert S\vert, k)/k1$ for every equivalence class 
 $S\in\ov{\alpha}/\approx$ of a class of parallel arrows $\ov{\alpha}$ such 
 that $\vert S\vert\ge2$ if $\car k\neq2$ and $\vert S\vert>2$ if $\car k=2$.
\end{Kor}

\begin{Kor}
 Let $Q$ be a connected quiver and $\Lambda=kQ/\langle Z\rangle$ a finite 
 dimensional monomial algebra such that $L_{-1}=k(Q_0//Q_1)\cap\Ker\psi_1=0$.
 The following conditions are equivalent:
 \begin{enumerate}
  \item The Lie algebra ${\rm H}^1(\Lambda,\Lambda)$ is solvable.
  \item Every equivalence class $S$ of a class of parallel arrows 
        $\ov{\alpha}\in Q_1//$ of $Q$ contains one and only one arrow if the 
        characteristic of the field $k$ is not 2. In the case $\car k=2$ we 
        have $\vert S\vert\le2$ for all $S\in\ov{\alpha}/\approx$, 
        $\ov{\alpha}\in Q_1//$.
 \end{enumerate}
\end{Kor}

Since in characteristic $0$ the Lie algebra  of a connected algebraic group is solvable if and only if the
connected algebraic group is solvable (as a group), using Proposition \ref{L-1} we infer the following
result also proved by Guil-Asensio and Saor\'{\i}n (see Corollary 2.22 in \cite{GAS2}):

\begin{Kor}
 Let $k$ be a field of characteristic $0$, $Q$ a connected quiver and $\Lambda=kQ/\langle Z\rangle$ a
 finite dimensional monomial algebra. Then the identity component of the algebraic group of outer
 automorphisms $\Out(\Lambda)^\circ$ is solvable if and only if every two parallel arrows of $Q$ are not
 equivalent.
\end{Kor}

\begin{Def}
  Let $\ov{Q}$ be the subquiver of $Q$ obtained by taking a representative 
 for every class of parallel arrows of $Q$.
\end{Def}

\begin{Thm}
\label{semisimple}
 Let $Q$ be a connected quiver and $\Lambda=kQ/\langle Z\rangle$ a finite dimensional monomial algebra
 such that $L_{-1}=k(Q_0//Q_1)\cap\Ker\psi_1=0$. The following conditions are equivalent:
 \begin{enumerate}
  \item The Lie algebra ${\rm H}^1(\Lambda,\Lambda)$ is semisimple.
  \item The quiver $\ov{Q}$ is a tree, $Q$ has at least one non trivial class of parallel arrows and the
        ideal $\langle Z\rangle$ is completely saturated. If the characteristic of the field $k$ is $2$,
    then $Q$ does not have a class of parallel arrows containing exactly two arrows.
  \item ${\rm H}^1(\Lambda,\Lambda)$ is isomorphic to the non trivial product of Lie algebras
        $$\prod_{\ov{\alpha}\in Q_1//}\mathfrak{pgl}(\vert\ov{\alpha}\vert,k)$$
    where $\vert\ov{\alpha}\vert\neq2$ if the characteristic of $k$ is equal to $2$.
 \end{enumerate}
\end{Thm}

\begin{Bew}
 Let $\mathcal{B}$ be a basis of $L_0$ as described in Proposition \ref{basis}.
 
 {\it (i)}\ $\Rightarrow$\ {\it (ii)}\ : If $\ov{Q}$ contained a (non oriented or oriented) cycle, there
 would exist a class $\ov{\alpha}\in Q_1//$ of parallel arrows such that 
 $(\alpha,\alpha)\in\mathcal{B}$ for all $\alpha\in\ov{\alpha}$ (see 
 \ref{full})
 and so $k(\sum_{\alpha\in\ov{\alpha}}(\alpha,\alpha))$ would be a nontrivial
 abelian ideal of $L_0$ which contradicts the fact that
 $0=\Rad{\rm H}^1(\Lambda,\Lambda)=\Rad L_0\oplus\bigoplus_{i\ge1}L_i$. 
 To insure that ${\rm H}^1(\Lambda,\Lambda)=L_0\neq0$ it is necessary 
 that $Q$ contains at least one nontrivial class of parallel
 arrows. Suppose that $\langle Z\rangle$ is not completely saturated. According to Remark \ref{saturated}
 there exist two parallel arrows which are not equivalent. This means that there
 exists a class of parallel arrows $\ov{\alpha}\in Q_1//$ such that $\ov{\alpha}/\approx$ contains at
 least two elements, say $R$ and $S$. The properties of the basis $\mathcal{B}$ show that either 
 $\sum_{a\in R}(a,a)\neq0$ or $\sum_{a\in S}(a,a)\neq0$ which contradicts 
 the fact that
 $\Rad L_0\subset\Rad{\rm H}^1(\Lambda,\Lambda)=0$. The last statement is 
 an immediate consequence of Theorem \ref{radical}. 
 
 {\it (ii)}\ $\Rightarrow$\ {\it (iii)}\ : Let $Q$ and $Z$ satisfy the given conditions. From the fact that 
 $\ov{Q}$ is a tree, we deduce $\bigoplus_{i\ge1}L_i=0$ and so ${\rm H}^1(\Lambda,\Lambda)=L_0$. Since
 $\langle Z\rangle$ is completely saturated we have 
 $L_0=k(Q_1//Q_1)/\langle\sum_{a\in Q_1e}(a,a)-\sum_{a\in eQ_1}(a,a)\mid 
   e\in Q_0\rangle$. This and the fact that $Q$ does not have an equivalence class containing two arrows if
   $\car k=2$, show that its radical is generated by the elements of type 
   $\sum_{\alpha\in\ov{\alpha}}(\alpha,\alpha)$, $\ov{\alpha}\in Q_1//$, (see Theorem \ref{radical}). These
   elements being equal to zero we see that $\Rad{\rm H}^1(\Lambda,\Lambda)=0$ and Corollary \ref{pgl}
   yields the result.
 
 {\it (iii)}\ $\Rightarrow$\ {\it (i)}\ : Clear in view of the following Remark.
\end{Bew}

\begin{Bem}
 The Lie algebra $\mathfrak{pgl}(n,k):=\mathfrak{gl}(n,k)/k1$, $n\ge2$, is isomorphic to the classical simple Lie
 algebra $\mathfrak{sl}(n,k)$ of $n\times n$-matrices having trace zero if the characteristic of the field $k$
 does not divide $n$. If $\car k$ divides $n$ and $n\neq2$, then $\mathfrak{pgl}(n,k)$ is a semisimple algebra
 without being a direct product of simple Lie algebras.
\end{Bem}

\begin{Kor}
\label{simple}
 Let $Q$ be a connected quiver and $\Lambda=kQ/\langle Z\rangle$ a finite dimensional monomial algebra
 such that $L_{-1}=k(Q_0//Q_1)\cap\Ker\psi_1=0$. The following conditions are equivalent:
 \begin{enumerate}
  \item The Lie algebra ${\rm H}^1(\Lambda,\Lambda)$ is simple.
  \item The quiver $\ov{Q}$ is a tree, the quiver $Q$ has exactly one class of parallel arrows
        $\ov{\alpha}=\{\alpha_1,\ldots,\alpha_n\}$ such that $n\ge2$ and the characteristic of the field $k$
    does not divide $n$. The ideal $\langle Z\rangle$ is completely saturated.  
    
  \item There exists an integer $n\ge2$ such that the characteristic of the field $k$ does not divide $n$ and
        ${\rm H}^1(\Lambda,\Lambda)$ is isomorphic to the Lie algebra $\mathfrak{sl}(n,k)$.
 \end{enumerate}
\end{Kor}

Since the first Hochschild cohomology group is the Lie algebra of the algebraic group of outer automorphisms
in characteristic $0$ (see Proposition \ref{Out}), we get the following Corollary (see Corollary 4.9 in \cite{GAS2} for a
partial result on hereditary algebras):

\begin{Kor}
 Let $k$ be a field of characteristic $0$, $Q$ a connected quiver and $\Lambda=kQ/\langle Z\rangle$ a finite
 dimensional monomial algebra. The following conditions are equivalent:
 \begin{enumerate}
  \item The algebraic group $\Out(\Lambda)^\circ$ is semisimple.
  \item The quiver $\ov{Q}$ is a tree, the quiver $Q$ has at least one non trivial class of parallel arrows
        and the ideal $\langle Z\rangle$ is completely saturated.
  \item The algebraic group $\Out(\Lambda)^\circ$ is isomorphic to the non trivial product of algebraic
        groups $\prod_{\ov{\alpha}\in Q_1//}{\bf PGL}(\vert\ov{\alpha}\vert,k)$.
 \end{enumerate}
\end{Kor}

\begin{Bew}
 In characteristic $0$ a connected algebraic group is semisimple if and only if 
 its Lie algebra is semisimple (see \cite{Hu} 13.5).
 
 {\it (i)}\ $\Leftrightarrow$\ {\it (ii)}\ : This equivalence follows immediately from the preceding Theorem
 and Proposition \ref{Out}.
 
 {\it (ii)}\ $\Rightarrow$\ {\it (iii)}\ : Let $Q$ and $Z$ satisfy the given condition. According to 
 \ref{semidirect} we know that $\Out(\Lambda)^\circ=\frac{H_\Lambda\cap U_\Lambda}{H_\Lambda\cap\Inn^\ast(\Lambda)}
 \rtimes\frac{V_\Lambda^\circ}{\hat{E}}$. Since $\ov{Q}$ is a tree, there is no path in $B':=B-Q_0\cup Q_1$ of
 length $\ge2$ which is parallel to an arrow. From the fact that we have for every element 
 $\sigma\in H_\Lambda\cap U_\Lambda$ that $\sigma(a)=a+\sum_{(a,\gamma)\in Q_1//B'}\lambda_{(a,\gamma)}\gamma$
 for all $a\in Q_1$ we deduce that $H_\Lambda\cap U_\Lambda=\{\id_\Lambda\}$. As the ideal $\langle Z\rangle$ is
 completely saturated by assumption, any two parallel arrows are equivalent. Therefore the main Theorem 2.20 of the
 article \cite{GAS2} implies $V_\Lambda^\circ\simeq\prod_{\ov{\alpha}\in Q_1//}{\bf GL}(\vert\ov{\alpha}\vert,k)$.
 Since the quiver $\ov{Q}$ is assumed to be a tree, we have 
 $\hat{E}={\rm Ch}(Q,k)=\prod_{\ov{\alpha}\in Q_1//}k^\ast I_{\vert\ov{\alpha}\vert}$ thanks to example 8
 in \cite{GAS1}. Hence we obtain
 $$\begin{array}{rcl}
   \Out(\Lambda)^\circ&=&\frac{H_\Lambda\cap U_\Lambda}{H_\Lambda\cap\Inn^\ast(\Lambda)}\rtimes
                      \frac{V_\Lambda^\circ}{\hat{E}}\ =\ \frac{V_\Lambda^\circ}{\hat{E}}\\
   &\simeq&\frac{\prod_{\ov{\alpha}\in Q_1//}{\bf GL}(\vert\ov{\alpha}\vert,k)}
              {\prod_{\ov{\alpha}\in Q_1//}k^\ast I_{\vert\ov{\alpha}\vert}}
   \ =\ \prod_{\ov{\alpha}\in Q_1//}\frac{{\bf GL}(\vert\ov{\alpha}\vert,k)}{k^\ast I_{\vert\ov{\alpha}\vert}}
   \ =\ \prod_{\ov{\alpha}\in Q_1//}{\bf PGL}(\vert\ov{\alpha}\vert,k)
 \end{array}$$

 {\it (iii)}\ $\Rightarrow$\ {\it (i)}\ : The semisimplicity of the Lie algebra $\mathcal{L}(\Out(\Lambda)^\circ)=
   \prod_{\ov{\alpha}\in Q_1//}\mathfrak{pgl}(\vert\ov{\alpha}\vert,k)$ implies the semisimplicity of the connected 
   algebraic group $\Out(\Lambda)^\circ$.
\end{Bew}

The fact that in characteristic $0$ the Lie algebra of a connected algebraic group is simple if and only if
the connected algebraic group is almost simple, yields the following result:

\begin{Kor}
 Let $k$ be a field of characteristic $0$, $Q$ a connected quiver and $\Lambda=kQ/\langle Z\rangle$ a 
 finite dimensional monomial algebra. The following conditions are equivalent:
 \begin{enumerate}
  \item The algebraic group $\Out(\Lambda)^\circ$ is almost simple.
  \item The quiver $\ov{Q}$ is a tree, the quiver $Q$ has exactly one class of parallel arrows
        $\ov{\alpha}=\{\alpha_1,\ldots\alpha_n\}$ of order $n\ge2$ and the ideal $\langle Z\rangle$ is
    completely saturated.
  \item There exists an integer $n\ge2$ such that the algebraic group $\Out(\Lambda)^\circ$ is isomorphic to
        the algebraic group ${\bf PGL}(n,k)$.
 \end{enumerate}
\end{Kor}

In order to get criteria for the commutativity and the reductivity of the 
Lie algebra ${\rm H}^1(\Lambda,\Lambda)$ we need to study its center. 
Therefore we introduce the following definitons:
for every class of parallel arrows $\ov{\alpha}$ of $Q$ we call a set $C\subset\ov{\alpha}$ connected, if
for every two arrows $\alpha_1$ and $\alpha_r$ of $C$ there exist arrows $\alpha_2,\ldots,\alpha_{r-1}\in C$ such
that we have $\alpha_i\le_{\langle Z\rangle}\alpha_{i+1}$ or $\alpha_{i+1}\le_{\langle Z\rangle}\alpha_i$
for all $i\in\{1,\ldots,r-1\}$. A connected set $C\subset\ov{\alpha}$ is called a connected component of
$\ov{\alpha}$ if it is maximal for the connection, \ie for every arrow $\beta\in\ov{\alpha}-C$ there is no
arrow $\alpha\in C$ such that $\alpha\le_{\langle Z\rangle}\beta$ or $\beta\ge_{\langle Z\rangle}\alpha$.
Clearly the connected components of a class of parallel arrows $\ov{\alpha}$ form a partition of
$\ov{\alpha}$.
 
\begin{Lem}
 Let $Q$ be a connected quiver and $\Lambda=kQ/\langle Z\rangle$ a finite dimensional 
 monomial algebra. 
 \begin{enumerate}
  \item The center ${\rm Z}\,(L_0)$ of the Lie algebra $L_0$ is generated by the elements 
   $\sum_{a\in C}(a,a)$ where $C$ denotes a connected component of a class of parallel arrows 
   of $Q$. 
  \item If the field $k$ has characteristic $0$ or if the quiver $Q$ does not have an
   oriented cycle, then the center ${\rm Z}\,({\rm H}^1(\Lambda,\Lambda))$ 
   of the Lie algebra ${\rm H}^1(\Lambda,\Lambda)$ is contained in the 
   center ${\rm Z}\,(L_0)$ of $L_0$.
 \end{enumerate}
\end{Lem}

\begin{Bew}
 {\it (i)}\ : For every element $(a,b)\in Q_1//Q_1$, $a\neq b$, we have 
 $[(a,b),(a,a)]=(a,b)$ so that ${\rm Z}\,(L_0)$ 
 is contained in the abelian Lie subalgebra of $L_0$ generated by the elements $(a,a)$, $a\in Q_1$. For
 every linear combination $\sum_{a\in Q_1}\lambda_a(a,a)$, $\lambda_a\in k$, we have 
 $[\sum_{a\in Q_1}\lambda_a(a,a),(b,c)]=(-\lambda_b+\lambda_c)(b,c)$ for all $(b,c)\in L_0$. Therefore
 $\sum_{a\in Q_1}\lambda_a(a,a)$ is contained in ${\rm Z}\,(L_0)$ if and only if we have 
 $\lambda_b=\lambda_c$ for
 all arrows $b,c$ such that $c\ge_{\langle Z\rangle}b$. This shows that ${\rm Z}\,(L_0)$ is generated by the
 elements $\sum_{a\in C}(a,a)$.
  
 {\it (ii)}\ :  In both cases we have $L_{-1}=0$ according to Proposition 
 \ref{L-1}. If the characteristic of $k$ is $0$, 
 then we have for every element 
 $0\neq\sum_{(a,\gamma)\in Q_1//B_{i+1}}\lambda_{(a,\gamma)}(a,\gamma)\in L_i$, 
 $\lambda_{(a,\gamma)}\in k$, $i\ge1$, that
 $$[\sum_{(a,\gamma)\in Q_1//B_{i+1}}\lambda_{(a,\gamma)}(a,\gamma),\sum_{b\in Q_1}(b,b)]=
   -i\sum_{(a,\gamma)\in Q_1//B_{i+1}}\lambda_{(a,\gamma)}(a,\gamma)\neq0$$
 This shows ${\rm Z}\,({\rm H}^1(\Lambda,\Lambda))\subset{\rm Z}\,(L_0)$.
 We assume now that $Q$ does not have an oriented cycle. Then every path $\gamma\in B$ of
 length $\ge2$ parallel to an arrow $a$ cannot contain $a$ and so 
 $[(a,\gamma),(a,a)]=(a,\gamma)$. Taking into account the fact
 ${\rm H}^1(\Lambda,\Lambda)=L_0\oplus\bigoplus_{i\ge1}k(Q_1//B_{i+1})\cap\Ker\psi_1$
 yields ${\rm Z}\,({\rm H}^1(\Lambda,\Lambda))\subset{\rm Z}\,(L_0)$. 
\end{Bew}
  
Recall that a Lie algebra is called reductive if its radical and its center are equal.

\begin{Prop}
 Let $Q$ be a connected quiver and $\Lambda=kQ/\langle Z\rangle$ a finite dimensional monomial algebra. If
 $Q$ does not have an oriented cycle or if the field $k$ has characteristic $0$, then the following
 conditions are equivalent:
 \begin{enumerate}
  \item The Lie algebra ${\rm H}^1(\Lambda,\Lambda)$ is reductive.
  \item The Lie ideal $\bigoplus_{i\ge1}L_i$ is equal to $0$ and the relation $\le_{\langle Z\rangle}$ is
        symmetric on every class of parallel arrows. If the characteristic of $k$ is $2$, then there 
    exists no equivalence class $S\in\ov{\alpha}/\approx$ containing exactly two arrows of a class of 
    parallel arrows $\ov{\alpha}$.
 \end{enumerate}
\end{Prop}

\begin{Bew}
 {\it (i)}\ $\Rightarrow$\ {\it (ii)}\ : Let be 
 $\Rad{\rm H}^1(\Lambda,\Lambda)={\rm Z}\,({\rm H}^1(\Lambda,\Lambda))$. If $Q$ does not have an oriented
 cycle or if $\car k=0$, then $L_{-1}=0$ and 
 $\Rad{\rm H}^1(\Lambda,\Lambda)=\Rad L_0\oplus\bigoplus_{i\ge1}L_i$. 
 The preceding Lemma shows that the center of ${\rm H}^1(\Lambda,\Lambda)$ is included
 in $L_0$ and thus $\Rad{\rm H}^1(\Lambda,\Lambda)\subset L_0$. This implies 
 $\bigoplus_{i\ge1}L_i=0$ and therefore ${\rm H}^1(\Lambda,\Lambda)=L_0$. 
 In view of Theorem \ref{radical} and the preceding Lemma, the assumption 
 $\Rad{L_0}=\Rad{\rm H}^1(\Lambda,\Lambda)={\rm Z}\,({\rm H}^1(\Lambda,\Lambda))={\rm Z}\,(L_0)$ 
 implies that there is no element $(a,b)\in L_0$ such that $(b,a)\not\in L_0$, \ie 
 the relation $\le_{\langle Z\rangle}$ is
 symmetric. Furthermore there exists no equivalence class $S\in\ov{\alpha}/\approx$, 
 $\ov{\alpha}\in Q_1//$, containing exactly two elements if $\car k=2$.
 
 {\it (ii)}\ $\Rightarrow$\ {\it (i)}\ : The assumptions imply $L_{-1}=0$. From 
 $\bigoplus_{i\ge1}L_i=0$ we deduce ${\rm H}^1(\Lambda,\Lambda)=L_0$. 
 Theorem \ref{radical} and the preceding Lemma show that 
 $\Rad{\rm H}^1(\Lambda,\Lambda)=\Rad L_0={\rm Z}\,(L_0)={\rm Z}\,({\rm H}^1(\Lambda,\Lambda))$.
\end{Bew}

Since in characteristic $0$ a connected algebraic group is reductive if and only if its Lie algebra is
reductive, we deduce immediately using Proposition \ref{Out} and 
Corollary \ref{dico}

\begin{Kor}
 Let $k$ be a field of characteristic $0$, $Q$ a connected quiver and 
 \linebreak
 $\Lambda=kQ/\langle Z\rangle$ a finite
 dimensional monomial algebra. Then the identity component of the algebraic group of the outer automorphisms
 $\Out(\Lambda)^\circ$ is reductive if and only if the relation $\le_{\langle Z\rangle}$ is symmetric on
 every class of parallel arrows of $Q$ and if 
 $\frac{H_\Lambda\cap U_\Lambda}{H_\Lambda\cap\Inn^\ast(\Lambda)}=0$.
\end{Kor}

Define a sequence of ideals of a Lie algebra $L$ by setting $\mathcal{C}^0L:=L(=L^{(0)})$,
$\mathcal{C}^1L:=[L,L](=L^{(1)})$, $\mathcal{C}^2L:=[L,\mathcal{C}^1L]$, $\ldots$,
$\mathcal{C}^nL:=[L,\mathcal{C}^{n-1}L]$. A Lie algebra $L$ is called nilpotent if there exists a 
nonnegative
integer $m$ such that $\mathcal{C}^mL=0$. The integer $m$ such that $\mathcal{C}^mL=0$ and
$\mathcal{C}^{m-1}L\neq0$ is called the nilindex of $L$. A Lie algebra $L$ is called filiform if it is
nilpotent of maximal nilindex $m$ that is $\dim\mathcal{C}^kL={m-k}$ for $1\le k\le m$. The fact
$L^{(n)}\subset\mathcal{C}^nL$ for all $n\in\N$ implies that nilpotent algebras are solvable.

\begin{Prop}
 Let $Q$ be a connected quiver and $\Lambda=kQ/\langle Z\rangle$ a finite dimensional monomial algebra. If
 $Q$ does not have an oriented cycle or if $\car k=0$, then the following conditions are equivalent:
 \begin{enumerate}
  \item The Lie algebra ${\rm H}^1(\Lambda,\Lambda)$ is abelian.
  \item The Lie algebra ${\rm H}^1(\Lambda,\Lambda)$ is filiform.
  \item The Lie algebra ${\rm H}^1(\Lambda,\Lambda)$ is nilpotent.
  \item $\bigoplus_{i\ge1}L_i=0$ and there exist no parallel arrows $a\neq b$ satisfying $a\le_{\langle Z\rangle}b$.
  \item The Lie algebra ${\rm H}^1(\Lambda,\Lambda)$ is generated by the elements $(a,a)$ of $L_0$. 
  \item The dimension of the Lie algebra ${\rm H}^1(\Lambda,\Lambda)$ equals the Euler characteristic
        $\vert Q_1\vert-\vert Q_0\vert-1$.
 \end{enumerate}
\end{Prop}

\begin{Bew}
 {\it (i)}\ $\Rightarrow$\ {\it (ii)}\ and
 {\it (ii)}\ $\Rightarrow$\ {\it (iii)}\ are obvious.
 
 {\it (iii)}\ $\Rightarrow$\ {\it (iv)}\ : If there were two parallel arrows $a\neq b$ such that 
 $a\le_{\langle Z\rangle}b$, we would have $(b,a)\in L_0$ and $[(a,a),(b,a)]=(b,a)$, and so
 $(b,a)\in\mathcal{C}^n{\rm H}^1(\Lambda,\Lambda)$ for all $n\in\N$, contrary to the assumption. If $Q$
 does not have an oriented cycle, then there exists for every element 
 $0\neq\sum_{(a,\gamma)\in Q_1//B_{i+1}}\lambda_{(a,\gamma)}(a,\gamma)\in L_i$, $i\ge1$, 
 an arrow $b\in Q_1$ such that
 $[\sum_{(a,\gamma)\in Q_1//B_{i+1}}\lambda_{(a,\gamma)}(a,\gamma),(b,b)]=
  \sum_{(b,\gamma)\in Q_1//B_{i+1}}\lambda_{(b,\gamma)}(b,\gamma)\neq0$ 
 and $[\sum_{(b,\gamma)\in Q_1//B_{i+1}}\lambda_{(b,\gamma)}(b,\gamma),(b,b)]=
 \sum_{(b,\gamma)\in Q_1//B_{i+1}}\lambda_{(b,\gamma)}(b,\gamma)\neq0$.
 Thus we have 
 $0\neq\sum_{(b,\gamma)\in Q_1//B_{i+1}}\lambda_{(b,\gamma)}(b,\gamma)
 \in\mathcal{C}^n{\rm H}^1(\Lambda,\Lambda)$
 for all $n\in\N$. Since the Lie algebra ${\rm H}^1(\Lambda,\Lambda)$ is 
 assumed to be nilpotent, it
 follows that $\bigoplus_{i\ge1}L_i=0$. In 
 case $\car k=0$ we have for every element 
 $\sum_{(a,\gamma)\in Q_1//B_{i+1}}\lambda_{(a,\gamma)}(a,\gamma)\in L_i$, $i\ge1$, 
 $$[\sum_{(a,\gamma)\in Q_1//B_{i+1}}\lambda_{(a,\gamma)}(a,\gamma),\sum_{b\in Q_1}(b,b)]=
   -i\sum_{(a,\gamma)\in Q_1//B_{i+1}}\lambda_{(a,\gamma)}(a,\gamma)\neq0$$
 and thus 
 $\sum_{(a,\gamma)\in Q_1//B_{i+1}}\lambda_{(a,\gamma)}(a,\gamma)\in\mathcal{C}^n{\rm H}^1(\Lambda,\Lambda)$
 for all $n\in\N$ which implies $\bigoplus_{i\ge1}L_i=0$.
 
 {\it (iv)}\ $\Rightarrow$\ {\it (v)}\ : By assumption ${\rm H}^1(\Lambda,\Lambda)=L_0$ 
 and $L_0$ does not contain an element $(a,b)$, $a\neq b$.
 
 {\it (v)}\ $\Rightarrow$\ {\it (vi)}\ : We deduce from the given condition 
 that
 $\dim\Ker\psi_1=\vert Q_1\vert$. Since the dimension of 
 $\im\psi_0=\langle\sum_{a\in Q_1e}(a,a)-\sum_{a\in eQ_1}(a,a)\mid e\in Q_0\rangle$ is equal to 
 $\vert Q_0\vert-1$ we obtain $\dim{\rm H}^1(\Lambda,\Lambda)=\dim\Ker\psi_1-\dim\im\psi_0=
 \vert Q_1\vert-\vert Q_0\vert +1$.
 
 {\it (vi)}\ $\Rightarrow$\ {\it (i)}\ : From $\dim\im\psi_0=\vert Q_0\vert-1$ and
 $\dim{\rm H}^1(\Lambda,\Lambda)=\vert Q_1\vert -\vert Q_0\vert+1$ it follows that 
 $\dim\Ker\psi_1=\vert Q_1\vert$ and thus $\Ker\psi_1$ is generated by the elements $(a,a)\in Q_1//Q_1$.
 Since the bracket of the Lie algebra $k(Q_1//B)$ is such that $[(a,a),(b,b)]=0$ for all $a,b\in Q_1$, we
 see that ${\rm H}^1(\Lambda,\Lambda)$ is an abelian Lie algebra.
\end{Bew}

The following Corollary is clear thanks to Proposition \ref{Out}, because in characteristic $0$ a connected algebraic
group is abelian (resp. nilpotent) if and only if its Lie algebra is abelian (resp. nilpotent) (see
\cite{Hu} 13.4 and 10.5). It generalizes Guil-Asensio and Saor\'{\i}n's criterion on the 
commutativity of the
algebraic group $V_\Lambda^\circ(=\varepsilon_\Lambda(G_\Lambda))$ (see Corollary 2.23 in \cite{GAS2}).

\begin{Kor}
 Let $k$ be a field of characteristic $0$ and $\Lambda=kQ/\langle Z\rangle$ a finite dimensional monomial
 algebra. The following conditions are equivalent:
 \begin{enumerate}
  \item The algebraic group $\Out(\Lambda)^\circ$ is abelian.
  \item The algebraic group $\Out(\Lambda)^\circ$ is nilpotent.
  \item[{\it (iii)}] For every path $\gamma\in B$ parallel to an arrow $a\neq\gamma$ there exists a path $p$ of
        $Z$ such that $p^{(a,\gamma)}\neq0$. 
 \end{enumerate}
\end{Kor}

\section{Application to group algebras where the group admits a normal cyclic
         Sylow $p$-subgroup}

The starting point for the following application is the paragraph
`Repr\'esentations modulaires des groupes finis' in Gabriel's article
\cite{Ga} (see also \cite{GR} p.75). Let $k$ be an algebraically closed field 
and $G$ a finite group. If the characteristic of $k$ does not divide the order 
of $G$, then Maschke's Theorem states that the group algebra $kG$ is
semisimple. Since $k$ is supposed to be algebraically closed, this is 
equivalent to saying that the group algebra $kG$ is separable and thus 
${\rm H}^1(kG,kG)=0$. So henceforth let $k$ be a field of 
characteristic $p$ dividing the order 
of $G$. Let $n=p^aq$ be the order of $G$ with $q$ prime to $p$. We are
interested in this section in the particular case where $G$ contains only one
cyclic Sylow $p$-subgroup $S$, necessarily normal in $G$. Then $S$ is a 
normal Hall subgroup of $G$ and Schur's splitting theorem (see \cite{Sc} 
9.3.6) states that there
exists a supplement $K$ of $S$ in $G$, \ie a subgroup such that $S\cap K=\{1\}$
and $SK=G$. It is unique up to a group isomorphism. If $\sigma$ is a generator
of $S$, the action of $K$ on $S$ by conjugation is given by a formula of the 
type
$$x\sigma x^{-1}=\sigma^{\chi(x)}$$
with $x\in K$ and $\chi(x)\in(\Z/p^a\Z)^\times$. The group $G$ is the 
semidirect product of $S$ and $K$. For every $kK$-module $N$ we denote by 
${_\chi N}$ the underlying $k$-vector space of $N$ equipped with a new action
$\ast$ of $K$ such that we have
$$x\ast m:=\chi(x)x\cdot m$$
for all $x\in K$ and $m\in N$. If the $kK$-module $N$ belongs to the 
isomorphism class $e$, we denote by $\chi e$ the isomorphism class of the 
$kK$-module ${_\chi N}$. That which enables us to study the Lie algebra 
${\rm H}^1(kG,kG)$ is the association of a quiver $Q$ to the 
group $G$ in the following way: the set $Q_0$ of vertices of $Q$ consists in the set of 
isomorphism classes of simple $kK$-modules and for every vertex $e$ we take an
arrow $e\rightarrow\chi e$. The quiver $Q$ is a disjoint union of crowns 
\ie of oriented cycles. Gabriel proved in \cite{Ga} that the category of
$kG$-modules is equivalent to the category of 
$kQ/\langle Q_{p^a}\rangle$-modules. In other words the $k$-algebras $kG$ and
$kQ/\langle Q_{p^a}\rangle$ are Morita equivalent. Note that in general they 
are not isomorphic. A necessary and sufficient condition for the existence of
an algebra isomorphism between $kG$ and $kQ/\langle Q_{p^a}\rangle$ is the 
commutativity of the group $K$. Since the Hochschild cohomology 
${\rm H}^*(\Lambda,\Lambda)$ is for every $k$-algebra $\Lambda$ Morita 
invariant as a Gerstenhaber algebra (see \cite{GS} p.143), we conclude that the Lie algebras ${\rm H}^1(kG,kG)$ and 
${\rm H}^1(kQ/\langle Q_{p^a}\rangle,kQ/\langle Q_{p^a}\rangle)$ are 
isomorphic. Furthermore 
${\rm H}^1(kQ/\langle Q_{p^a}\rangle,kQ/\langle Q_{p^a}\rangle)$ is a product 
of Lie algebras 
${\rm H}^1(kC/\langle C_{p^a}\rangle,kC/\langle C_{p^a}\rangle)$ where $C$ 
is a crown of $Q$, because the quiver $Q$ is a disjoint union of crowns and the
Hochschild cohomology is additive on a product of algebras. Thus we have
$${\rm H}^1(kG,kG)=\prod_{C\ {\rm crown\ of}\ Q}
  {\rm H}^1(kC/\langle C_{p^a}\rangle,kC/\langle C_{p^a}\rangle)$$

\begin{Thm}
 Let $k$ be an algebraically closed field and $G$ a finite group such that
 the characteristic $p>0$ of $k$ divides the order of $G$. Let $n=p^aq$ be the
 order of $G$ with $q$ prime to $p$. Suppose that $G$ contains only one Sylow 
 $p$-subgroup which in addition is cyclic. The following conditions are 
 equivalent:
 \begin{enumerate}
  \item The Lie algebra ${\rm H}^1(kG,kG)$ is semisimple.
  \item The group $G$ is the direct product of a group $K$ of order $q$ and of 
        a cyclic group $C_p$ of order $p$. The characteristic $p$ of the field
        $k$ is different from $2$.
  \item ${\rm H}^1(kG,kG)$ is a product of Witt algebras 
        $W(1,1):=\Der(k[X]/(X^p))$ and the characteristic $p$ of the field $k$ 
        is different from $2$.
 \end{enumerate}
\end{Thm}

\begin{Bew}
 Let $Q$ be the above-described quiver associated to the group $G$.

 {\it (i)}\ $\Rightarrow$\ {\it (ii)}\ : It is clear that the Lie algebra 
 ${\rm H}^1(kG,kG)$ is semisimple if and only if the Lie algebra
 ${\rm H}^1(kC/\langle C_{p^a}\rangle,kC/\langle C_{p^a}\rangle)$ is semisimple
 for every crown $C$ of the quiver $Q$. From Proposition \ref{L-1} and 
 Theorem \ref{semisimple} we deduce that $Q$ does not have a crown of 
 length $\ge2$, \ie $Q$ is a disjoint union of loops. Therefore we have 
 $\chi(x)=1$ for all $x\in G$ which implies that $G$ is the direct product of
 the Sylow $p$-subgroup $S$ and its supplement $K$. Proposition  
 \ref{loop} shows that for every loop $C$ of the quiver $Q$ the Lie algebra 
 ${\rm H}(kC/\langle C_{p^a}\rangle,kC/\langle C_{p^a}\rangle)$ is semisimple 
 if and only if $p^a=p>2$.

 {\it (ii)}\ $\Rightarrow$\ {\it (iii)}\ : Let $G=K\times C_p$ and so $\chi=1$.
 Therefore $Q$ is a disjoint union of loops. According to Proposition
 \ref{loop} we have for every loop $C$ of the quiver $Q$ that
 ${\rm H}^1(kC/\langle C_{p^a}\rangle,kC/\langle C_{p^a}\rangle)$ is isomorphic
 to the Witt Lie algebra $W(1,1)$.

 {\it (iii)}\ $\Rightarrow$\ {\it (i)}\ : This is clear, because the Witt 
 algebra is  one of the nonclassical simple Lie algebras if the characteristic 
 $p>0$ of $k$ is different from $2$.
\end{Bew}


\begin{Kor}
 Let $k$ be an algebraically closed field and $G$ a finite group such that
 the characteristic $p>0$ of $k$ divides the order of $G$. Let $n=p^aq$ be the
 order of $G$ with $q$ prime to $p$. Suppose that $G$ contains only one Sylow 
 $p$-subgroup which in addition is cyclic. The following conditions are 
 equivalent:
 \begin{enumerate}
  \item The Lie algebra ${\rm H}^1(kG,kG)$ is simple.
  \item The group $G$ is cyclic of order $p$ and the characteristic $p$ of the
        field $k$ is different from $2$.
  \item ${\rm H}^1(kG,kG)$ is isomorphic to the Witt Lie algebra
        $W(1,1):=\Der(k[X]/(X^p))$ and the characteristic $p$ of the field $k$
        is different from $2$.
 \end{enumerate}
\end{Kor}

\begin{Bew}
 Let $Q$ be the above-described quiver associated to the group $G$.

 {\it (i)}\ $\Rightarrow$\ {\it (ii)}\ : Let the Lie algebra ${\rm H}^1(kG,kG)$
 be simple and thus semisimple. The preceding Theorem shows that $G$ is the 
 direct product of a group $K$ of order $q$ and of a cyclic group $C_p$ of 
 order $p$, where $q$ and $p$ are prime. Therefore we have $\chi=1$ which
 implies that the quiver $Q$ is a disjoint union of loops. From the fact that
 ${\rm H}^1(kG,kG)=
  \prod{\rm H}^1(kC/\langle C_p\rangle, kC/\langle C_p\rangle)$ is a simple 
 Lie algebra it follows that $Q$ has only one loop. Since the field $k$ is
 algebraically closed and since $p$ does not divide the order $q$ of the group 
 $K$, the number of isomorphism classes of simple $kK$-modules, which is the
 number of vertices of $Q$, is equal to the number of
 conjugation classes of the group $K$. Thanks to the fact that the 
 loop $Q$ has only one vertex we obtain that $K$ has only one conjugation 
 class and so $K=\{1\}$.  

 {\it (ii)}\ $\Rightarrow$\ {\it (iii)}\ : Since the quiver $Q$ associated to 
 the cyclic group $C_p$ is the loop, Proposition \ref{loop} shows that
 ${\rm H}^1(kG,kG)$ is isomorphic to the Witt Lie algebra.

 {\it (iii)}\ $\Rightarrow$\ {\it (i)}\ : clear
\end{Bew}

\end {document}